\documentclass[11pt]{article}
\usepackage[usenames]{color}
\usepackage{amssymb,psfig,epsfig,latexsym,graphicx}

\usepackage{amssymb}
\usepackage[colorlinks=true,
linkcolor=webgreen,
filecolor=webbrown,
citecolor=webgreen]{hyperref}

\definecolor{webgreen}{rgb}{0,.5,0}
\definecolor{webbrown}{rgb}{.6,0,0}

\newtheorem{theorem}{Theorem}
\newtheorem{lemma}[theorem]{Lemma}
\newtheorem{coro}[theorem]{Corollary}

\newtheorem{obs}[theorem]{Observation}

\newcommand{\eqn}[1]{(\ref{#1})}
\newcommand{\bsq}{{\vrule height .9ex width .8ex depth -.1ex }}

\newcommand{\NN}{{\mathbb N}}

\newcommand{\PP}{{\mathbb P}}

\newcommand{\sC}{{\mathcal C}}

\newcommand{\eeq}{\end{equation}}
\newcommand{\Snj}{\scriptstyle}
\newcommand{\beql}[1]{\begin{equation}\label{#1}}

\newcommand{\bone}{\mbox{\boldmath $1$}}
\newcommand{\btwo}{\mbox{\boldmath $2$}}
\newcommand{\bthree}{\mbox{\boldmath $3$}}
\newcommand{\bfour}{\mbox{\boldmath $4$}}
\newcommand{\bfive}{\mbox{\boldmath $5$}}
\newcommand{\bsix}{\mbox{\boldmath $6$}}

\newcommand{\beight}{\mbox{\boldmath $8$}}
\newcommand{\bnine}{\mbox{\boldmath $9$}}
\newcommand{\bten}{\mbox{\boldmath $10$}}
\newcommand{\beleven}{\mbox{\boldmath $11$}}
\newcommand{\btf}{\mbox{\boldmath $24$}}
\newcommand{\btt}{\mbox{\boldmath $32$}}
\newcommand{\bss}{\mbox{\boldmath $67$}}
\newcommand{\bonf}{\mbox{\boldmath $195$}}
\newcommand{\bons}{\mbox{\boldmath $196$}}
\newcommand{\bfeh}{\mbox{\boldmath $580$}}
\newcommand{\bfeo}{\mbox{\boldmath $581$}}

\newcommand{\ubtwo}{\mbox{\underline{\boldmath $2$}}}
\newcommand{\ubthree}{\mbox{\underline{\boldmath $3$}}}
\newcommand{\ubfour}{\mbox{\underline{\boldmath $4$}}}

\newcommand{\cto}{\mbox{\textcircled{2}}{}_0}

\makeatletter
\def\@sect#1#2#3#4#5#6[#7]#8{\ifnum #2>\c@secnumdepth
     \def\@svsec{}\else
     \refstepcounter{#1}\edef\@svsec{\csname the#1\endcsname.\hskip .75em }\fi
     \@tempskipa #5\relax
      \ifdim \@tempskipa>\z@
        \begingroup #6\relax
          \@hangfrom{\hskip #3\relax\@svsec}{\interlinepenalty \@M #8\par}%
        \endgroup
       \csname #1mark\endcsname{#7}\addcontentsline
         {toc}{#1}{\ifnum #2>\c@secnumdepth \else
                      \protect\numberline{\csname the#1\endcsname}\fi
                    #7}\else
        \def\@svsechd{#6\hskip #3\@svsec #8\csname #1mark\endcsname
                      {#7}\addcontentsline
                           {toc}{#1}{\ifnum #2>\c@secnumdepth \else
                             \protect\numberline{\csname the#1\endcsname}\fi
                       #7}}\fi
     \@xsect{#5}}
\def\@begintheorem#1#2{\it \trivlist \item[\hskip \labelsep{\bf #1\ #2.}]}

\def\section{\@startsection {section}{1}{\z@}{-3.5ex plus -1ex minus
 -.2ex}{2.3ex plus .2ex}{\normalsize\bf}}
\makeatother
\makeatletter
\def\subsection{\@startsection {subsection}{1}{\z@}{-3.5ex plus -1ex minus
 -.2ex}{2.3ex plus .2ex}{\normalsize\bf}}

\makeatother

\def\am{a^{(m)}}
\def\Am{A^{(m)}}
\def\bm{\beta^{(m)}}
\def\sm{\sigma^{(m)}}
\def\tm{\tau^{(m)}}
\def\Bnm{B_n^{(m)}}

\begin{document}


\begin{center}
{\large\bf A Slow-Growing Sequence Defined by an Unusual Recurrence} \\
\vspace*{+.2in}

Fokko J. van de Bult${}^{(a)}$,
Dion C. Gijswijt${}^{(a)}$,
John P. Linderman${}^{(b)}$, \\
N. J. A. Sloane${}^{(b)}$ and
Allan R. Wilks${}^{(b)}$ \\
\vspace*{+.2in}
${}^{(a)}$Korteweg-de Vries Institute for Mathematics,
University of Amsterdam, \\
Plantage Muidergracht 24,
1018 TV Amsterdam, Netherlands. \\
${}^{(b)}$AT\&T Shannon Labs,
180 Park Avenue,
Florham Park, NJ 07932--0971, USA.
\vspace*{+.2in}

Email: fjvdbult@science.uva.nl,
dion.gijswijt@gmail.com,
jpl@research.att.com,
njas@research.att.com,
allan@research.att.com. \\
\vspace*{+.2in}
June 26, 2004; last revised February 22, 2006
\vspace*{+.2in}


{\bf Abstract}
\end{center}

The sequence starts with $a(1) =1$; to extend it one writes the sequence
so far as $XY^k$, where $X$ and $Y$ are strings of integers, $Y$ is
nonempty and $k$ is as large as possible: then the next term is $k$.
The sequence begins 1, 1, 2, 1, 1, 2, 2, 2, 3, 1, 1, 2, 1, 1, 2, 2, 2, 3,
2, $\ldots$ A $4$ appears for the first time at position 220, but a $5$
does not appear until about position $10^{10^{23}}$.  
The main result of the paper is a proof that the sequence is unbounded.
We also present results from extensive numerical investigations
of the sequence and of certain derived sequences, culminating
with a heuristic argument that $t$
(for $t=5,6, \ldots$) appears for the first time at about position 
$2\uparrow (2\uparrow (3\uparrow (4\uparrow (5\uparrow \ldots 
\uparrow ({(t-2)}\uparrow {(t-1)})))))$,
where $\uparrow$ denotes exponentiation.
The final section discusses generalizations.


\section{Introduction}\label{S1}
This paper introduces an integer sequence $A= a(1)$, $a(2)$, $a(3) ,
\ldots$ with some remarkable properties.  Define the {\em curling number}
$\sC (U)$ of a string $U = u(1)$, $u(2), \ldots , u(n)$ over some alphabet
$\Omega$ to be the largest integer $k \ge 1$ such that
\beql{EqM1}
U = X \underbrace{YY \ldots Y}_{\mbox{$k$ copies}} = XY^k ~,
\eeq
where $X$ and $Y$ are strings over $\Omega$ and $Y$ is nonempty.
Our sequence is defined by
\beql{EqM2}
a(1) =1, ~ a(n+1) = \sC (a(1), \ldots , a(n)) \quad\mbox{for} \quad
n \ge 1 \,.
\eeq

Then $a(2) = \sC (1) =1$, since we can only take $X$ to be the empty
string $\emptyset$, $Y=1$ and $k=1$; $a(3) = \sC (1,1) =2$, by taking $X =
\emptyset$, $Y=1$, $k=2$; $a(4) = \sC (1,1,2) =1$, by taking $X =1,1$,
$Y=2$, $k=1$ (as this example shows, there may be more than one choice
for $Y$); and so on.  The first 220 terms of $A$ are shown 
in Tables~\ref{T1} and \ref{T2}.

To avoid any possible confusion, for example with the ``Say What You See''
sequence studied in \cite{Con87}, we emphasize that the curling number
does {\em not} depend on the decimal representation of its arguments.
For example, if $U = (8,9,10,11,11,11)$, $\sC (U) =3$.

\begin{table}[htb]
$$
\begin{array}{cccccccccccccccccccc}
1 & 1 & \ubtwo \\
1 & 1 & \btwo & \ubtwo & \ubtwo & \ubthree \\
1 & 1 & 2 \\
1 & 1 & \btwo & \btwo & \btwo & \bthree & \ubtwo \\
1 & 1 & 2 \\
1 & 1 & 2 & 2 & 2 & 3 \\
1 & 1 & 2 \\
1 & 1 & \btwo & \btwo & \btwo & \bthree & \btwo & \ubtwo & \ubtwo & \ubthree & \ubtwo & \ubtwo & \ubtwo & \ubthree & \ubthree & \ubtwo \\
1 & 1 & 2 \\
1 & 1 & 2 & 2 & 2 & 3 \\
1 & 1 & 2 \\
1 & 1 & 2 & 2 & 2 & 3 & 2 \\
1 & 1 & 2 \\
1 & 1 & 2 & 2 & 2 & 3 \\
1 & 1 & 2 \\
1 & 1 & \btwo & \btwo & \btwo & \bthree & \btwo & \btwo & \btwo & \bthree & \btwo & \btwo & \btwo & \bthree & \bthree & \btwo & \ubtwo & \ubtwo & \ubthree & \ubtwo \\
\end{array}
$$

\caption{The first $98$ terms of the sequence.
In the notation to be introduced in Section \ref{S2},
the five underlined strings are the glue strings
$S_1^{(1)}, S_2^{(1)}, \ldots, S_5^{(1)}$
and the five bold-face strings are $T_2^{(1)}, T_3^{(1)}, \ldots, T_6^{(1)}$. }
\label{T1}
\end{table}
\begin{table}[htb]
$$
\begin{array}{cccccccccccccccccccccccccccc}
1 & 1 & 2 \\
1 & 1 & 2 & 2 & 2 & 3 \\
1 & 1 & 2 \\
1 & 1 & 2 & 2 & 2 & 3 & 2 \\
1 & 1 & 2 \\
1 & 1 & 2 & 2 & 2 & 3 \\
1 & 1 & 2 \\
1 & 1 & 2 & 2 & 2 & 3 & 2 & 2 & 2 & 3 & 2 & 2 & 2 & 3 & 3 & 2 \\
1 & 1 & 2 \\
1 & 1 & 2 & 2 & 2 & 3 \\
1 & 1 & 2 \\
1 & 1 & 2 & 2 & 2 & 3 & 2 \\
1 & 1 & 2 \\
1 & 1 & 2 & 2 & 2 & 3 \\
1 & 1 & 2 \\
1 & 1 &
\btwo & \btwo & \btwo & \bthree &
\btwo & \btwo & \btwo & \bthree &
\btwo & \btwo & \btwo & \bthree &
\bthree & \btwo & \btwo & \btwo & \bthree &
\btwo & \ubtwo & \ubtwo & \ubthree &
\ubtwo & \ubtwo & \ubtwo & \ubthree & \ubthree \\
~ & ~ & \ubtwo & \ubtwo & \ubtwo & \ubthree & \ubtwo & \ubtwo & \ubtwo & \ubthree & \ubtwo & \ubtwo & \ubtwo & \ubthree & \ubthree & \ubthree & \ubthree & \ubfour \\
\end{array}
$$

\caption{Terms 99 through 220 of the sequence,
up to the point where the first $4$ appears
($S_6^{(1)}$ is underlined, $T_7^{(1)}$ is shown in bold-face). }
\label{T2}
\end{table}

In Section \ref{S2} we describe the recursive structure of the sequence,
in particular explaining the block structure visible in Tables~\ref{T1}
and \ref{T2}.  The proof that this structure is valid is postponed to
Section \ref{S3}, where we give the main results of the paper, Theorems
\ref{lemF2} and \ref{Th1}.  Corollary \ref{CorFinite} shows that the sequence is unbounded.

In Section \ref{S4} we give empirical estimates for the lengths of the
blocks in the recursive structure, culminating in the estimate that $t \ge
5$ appears in the sequence for the first time at about position
\beql{Eq2}
2^{2^{\Snj 3^{\Snj 4^{\cdot^{\cdot^{\cdot^{\Snj {t-1}}}}}}}} \, ,
\eeq
a tower of height $t-1$.
These estimates are based on examination of the first two million terms
of the sequence $A$ and of the higher-order sequences $A^{(2)}$, 
$A^{(3)}$ and $A^{(4)}$ introduced in Section \ref{S2}.

The final section is devoted to comments and generalizations.
\S\ref{Sec51} discusses a certain plausible ``Finiteness Conjecture''
that arises from studying curling numbers.  
\S\ref{Sec52} discusses sequences that are obtained 
when the ``curling number transform'' (defined below)
is applied to certain well-known sequences.
Finally, \S\ref{Sec53} briefly mentions some
generalizations of our sequence, including a broad class of extensions
suggested by J. Taylor \cite{Tay04}.

Although the sequence $A$ grows very slowly, there are certainly familiar
sequences with an even slower growth rate, such as the inverse Ackermann
function \cite{Ack28}, the Davenport-Schinzel sequences \cite{ShAg95},
or the inverse to Harvey Friedman's sequence \cite{Fri01}.  Nevertheless,
we think the combination of slow growth, an unusual definition, and a
remarkable recursive structure makes the sequence noteworthy.

The sequence was invented by one of us (D.C.G.) while composing problems
for the Dutch magazine {\em Pythagoras}.  It now appears as sequence
A90822 in \cite{OEIS}.


\subsection*{Notation}
If $\Omega$ is a set, $\Omega^n$ denotes the strings of length $n$ from $\Omega$,
$\Omega^+$ is the set of all nonempty finite strings from $\Omega$,
and $\Omega^\ast$ is the set of all finite or infinite strings from $\Omega$,
including the empty string $\emptyset$.
Strings will usually be denoted by uppercase letters.
The elements of a string may or may not be separated by commas,
and a string may or may not be enclosed in parentheses.  
A sequence is an infinite string.
The length of $U \in \Omega^\ast$ (which may be $\infty$) will be denoted by $l(U)$.  

Products in $\Omega^\ast$ represent concatenation: 
if $U \in \Omega^+$, $V \in \Omega^\ast$ then $UV$ means $U$ followed by $V$.  
We will usually not concatenate two infinite strings.
A string
$U = u(1), \ldots, u(i)$ is said to be a {\em substring} of
$V = v(1), \ldots, v(j)$ if there is an $r$, $0 \le r \le j-i$,
such that $u(k) = v(k+r)$ for $k=1, \ldots, i$;
that is, if the elements of $U$ occur consecutively in $V$.
We say $V$ {\em contains} $U$ to indicate that $U$ is a substring of $V$.
Terms such as {\em prefix}, {\em suffix}, etc., have their
usual meanings --- see \cite{AlSh03} for formal definitions.  
A sequence $U$
is said to be a {\em subsequence}
of a sequence $V$ if $U$ can be obtained by deleting terms from $V$.

Usually $\Omega$ will be either
the nonnegative integers $\NN = \{0,1,2,3, \ldots \}$,
the positive integers $\PP = \{1,2,3, \ldots \}$,
or the set $\PP_m = \{m, m+1, m+2, \ldots \}$ for some integer $m \ge 1$.  

Given a sequence $U = u(1), u(2), \ldots \in \Omega^\ast$,
its {\em curling number transform} is the sequence
$U_\ast = u_\ast(1), u_\ast(2), \ldots \in \PP^\ast$
given by $u_\ast(1) = 1$ and
\beql{EqS3}
u_\ast(i) = \sC (u(1), \ldots, u(i-1)) \quad\mbox{for}\quad i \ge 2 \,.
\eeq
It is immediate from the definition \eqn{EqM2} that
our sequence $A$ is equal to its curling number transform,
and in fact is the unique sequence with this property.


\section{The recursive structure}\label{S2}
We introduce the notation in three stages: informally,
more formally and---in Section \ref{S3}---with a somewhat
different emphasis that will be needed to prove the
main theorems.

Informally, the sequence $A$ is built up recursively from
``blocks'' $B_n^{(1)}$ that are doubled at each step and are joined together by
``glue'' strings $S_n^{(1)}$.
When the glue strings alone are concatenated together they form
a sequence $A^{(2)}$ which has a similar structure to $A$: it is built up
recursively from blocks $B_n^{(2)}$ that are repeated three times at each step and
are joined together by ``second-order glue'' 
strings $S_n^{(2)}$.  When the second-order
glue strings are concatenated together they form a sequence $A^{(3)}$
which in turn has a similar structure, but now the blocks $B_n^{(3)}$ are repeated
four-fold at each step; and so on.  The proof that this description is
correct will be given in the next section.

We now make this description more precise.
The following description is correct, and is
the best way to think about the sequence.
However, we will not know for certain that it is correct
until the end of Section \ref{S3}.

The sequence $A$ is constructed from strings $B_n^{(1)}$ and $S_n^{(1)}$,
$n \ge 1$, which we call ``blocks'' and ``glue,'' respectively.
The initial block is $B_1^{(1)} =1$; the second block is $B_2^{(1)} =
B_1^{(1)} B_1^{(1)} S_1^{(1)} = 1~1~2$, where $S_1^{(1)} =2$; the third
block is
\begin{eqnarray*}
B_3^{(1)} & = & B_2^{(1)} B_2^{(1)} S_2^{(1)} \\
          & = & 1~1~2~1~1~2~2~2~3 \,,
\end{eqnarray*}
where $S_2^{(1)} = 2~2~3$, and so on, the $n$-th block for $n \ge 2$ being
\beql{EqR2}
B_n^{(1)} = B_{n-1}^{(1)} B_{n-1}^{(1)} S_{n-1}^{(1)} \,,
\eeq
where $S_{n-1}^{(1)}$ contains no 1's.
Then for all $n \ge 1$, $A$ begins with $B_n^{(1)}$ (and hence $A =
\displaystyle\lim_{n \to \infty} B_n^{(1)}$).

That is, for all $n \ge 2$, $A$ begins with two copies of $B_{n-1}^{(1)}$
followed by a ``glue'' string $S_{n-1}^{(1)}$ that contains no 1's.
$S_{n-1}^{(1)}$ is terminated by the first 1 that follows the initial
$B_{n-1}^{(1)} B_{n-1}^{(1)}$.  
Table~\ref{T1} shows $B_1^{(1)}$
through $B_6^{(1)}$ (the first row is $B_2^{(1)}$, the first two rows
together form
$B_3^{(1)}, \ldots$, and the whole table forms $B_6^{(1)}$), and
Tables~\ref{T1} and \ref{T2} together form $B_7^{(1)}$.  
The glue strings $S_1^{(1)}, S_2^{(1)}, \ldots, S_6^{(1)}$ are underlined.
By iterating \eqn{EqR2} we see that $B_n^{(1)}$ can also be written as
\beql{EqR2a}
B_n^{(1)} = 
B_{n-1}^{(1)} B_{n-2}^{(1)} \ldots
B_1^{(1)} B_1^{(1)} S_1^{(1)} S_2^{(1)} \ldots S_{n-1}^{(1)} \,.
\eeq
The terminating string $S_1^{(1)} S_2^{(1)} \ldots S_{n-1}^{(1)}$
(denoted by $T_n^{(1)}$ in Section \ref{S3})
is shown in bold-face in Tables~\ref{T1} and \ref{T2} for
$n=2, \ldots, 7$.

In Section \ref{S4} we state some conjectures about the lengths
of the blocks $B_n^{(1)}$ and of the glue strings $S_{n-1}^{(1)}$.
Assuming these conjectures are correct, $l(S_{n-1}^{(1)} )$ is much less than
$l(B_{n-1}^{(1)})$, and consequently $l(B_n^{(1)} )$ is roughly twice
$l(B_{n-1}^{(1)})$.

\begin{table}[ht]
$$
\begin{array}{ccccccccc}
2 & 2 & 2 & \ubthree \\
2 & 2 & 2 & 3 \\
2 & 2 & 2 & \bthree & \ubthree \\
2 & 2 & 2 & 3 \\
2 & 2 & 2 & 3 \\
2 & 2 & 2 & 3 & 3 \\
2 & 2 & 2 & 3 \\
2 & 2 & 2 & 3 \\
2 & 2 & 2 & \bthree & \bthree & \ubthree & \ubthree & \ubfour \\
2 & 2 & 2 & 3 \\
2 & 2 & 2 & 3 \\
2 & 2 & 2 & 3 & 3 \\
2 & 2 & 2 & 3 \\
2 & 2 & 2 & 3 \\
2 & 2 & 2 & 3 & 3 \\
2 & 2 & 2 & 3 \\
2 & 2 & 2 & 3 \\
2 & 2 & 2 & 3 & 3 & 3 & 3 & 4 \\
2 & 2 & 2 & 3 \\
2 & 2 & 2 & 3 \\
2 & 2 & 2 & 3 & 3 \\
2 & 2 & 2 & 3 \\
2 & 2 & 2 & 3 \\
2 & 2 & 2 & 3 & 3 \\
2 & 2 & 2 & 3 \\
2 & 2 & 2 & 3 \\
2 & 2 & 2 & \bthree & \bthree & \bthree & \bthree & \bfour & \ubthree \\
\end{array}
$$
\caption{The first $127$ terms of the second-order sequence $A^{(2)}$
(the successive second-order glue strings $S_1^{(2)}, S_2^{(2)}, S_3^{(2)}, S_4^{(2)}$ are underlined;
the strings  $T_2^{(2)}, T_3^{(2)}, T_4^{(2)}, T_5^{(2)}$ are shown in bold-face). }
\label{T3}
\end{table}

The above decomposition reduces the study of $A$ to the study of the
glue strings $S_n^{(1)}$.
We define the ``second-order sequence'' $A^{(2)} = a^{(2)} (1)$, $a^{(2)}
(2)$, $a^{(2)} (3), \ldots$ to be the concatenation $S_1^{(1)} S_2^{(1)}
S_3^{(1)} \ldots \in \PP_2^\ast$ of the glue strings.  
It will be shown later that $A^{(2)}$ can also be defined by
\begin{eqnarray}\label{EqM3}
&& a^{(2)} (1) =2 \,, \nonumber \\
&& a^{(2)} (n+1)
= \sC^{(2)} (a^{(2)} (1), a^{(2)} (2), \ldots, a^{(2)} (n))
\quad\mbox{for}\quad n \ge 1 \,,
\end{eqnarray}
where we define
\beql{EqM4}
\sC^{(m)} (U) = \max \{ m, \sC(U)\}
\eeq
for $m \ge 1$.  That is, if $\sC(U) =k$ is less than $m$ it is
``promoted'' to $m$
(we will say more about ``promotion'' at the end of Section \ref{S3}).
Of course $\sC^{(1)} = \sC$.

The first $127$ terms of $A^{(2)}$ are shown in Table~\ref{T3}, and the
reader can verify that they may indeed be obtained by starting with 2
and repeatedly applying the map $\sC^{(2)}$.

It is remarkable that $A^{(2)}$ has a similar structure to $A$, only now
the blocks are repeated three times.  That is, if we define $B_1^{(2)}=2$,
then for all $n \ge 2$, $A^{(2)}$ begins with a block
\beql{EqR3}
B_n^{(2)} = B_{n-1}^{(2)} B_{n-1}^{(2)} B_{n-1}^{(2)} S_{n-1}^{(2)} \,,
\eeq
consisting of three copies of $B_{n-1}^{(2)}$ followed by a ``second-order
glue'' string $S_{n-1}^{(2)} \in \PP_3^\ast$ that contains no 1's or 2's
and is terminated by the first number less than $3$ that follows the initial
$B_{n-1}^{(2)} B_{n-1}^{(2)} B_{n-1}^{(2)}$.  Table~\ref{T3} shows
$B_5^{(2)}$ (as well as $B_1^{(2)}$ through $B_4^{(2)}$).  
The glue strings $S_1^{(2)}, S_2^{(2)}, S_3^{(2)}, S_4^{(2)}$ are underlined.
$B_n^{(2)}$ ends with the string
$S_1^{(2)} S_2^{(2)} \ldots S_{n-1}^{(2)}$
(denoted by $T_n^{(2)}$ in Section \ref{S3});
these strings are shown in bold-face in Table~\ref{T3} for
$n=2, \ldots, 5$.

Again we have a conjectured estimate (see Section \ref{S4}) for the
lengths of the glue, which implies that $l(B_n^{(2)})$ is roughly three
times $l(B_{n-1}^{(2)} )$.

This analysis reduces the study of $A^{(2)}$ to the study of the
second-order glue strings $S_n^{(2)}$, and these, when concatenated, form
the third-order sequence $A^{(3)}$, which in turn has a similar structure.
And so on!


\section{The main theorems}\label{S3}
In this section we establish our main theorems, which
will show that the description
of the sequence given in Section \ref{S2} is correct.
To do this we must introduce our notation very carefully.
The following definitions (temporarily) supersede those in Section \ref{S2}.
For $m \ge 1$, the $m$th-order sequence
$\Am = a^{(m)}(1), a^{(m)}(2), a^{(m)}(3), \ldots \in \PP_m^{\ast}$ is defined by
\begin{eqnarray}\label{defA}
\am(1) & = & m \, , \nonumber \\
\am(i+1) & = & \sC^{(m)}(\am(1), \ldots, \am(i))
\quad\mbox{for} \quad i \ge 1 \,,
\end{eqnarray}
where $\sC^{(m)}$ is defined in \eqn{EqM4}.
Note that $A^{(1)}$ is our sequence $A$.
Theorem \ref{lemF2} will show
that $A^{(m+1)}$ is the concatenation of the glue strings for $\Am$.

For $m \ge 1$, $n \ge 1$,
the blocks $B_n^{(m)} \in \PP_m^\ast$
and the glue strings $S_n^{(m)} \in \PP_{m+1}^\ast$
are defined recursively, and independently of the $A^{(m)}$.
Corollary \ref{CorFinite} will show that all the strings $B_n^{(m)}$, $S_n^{(m)}$
and $T_n^{(m)}$ (defined below) are in fact finite, but at this point
we do not know that, and the definitions must allow for the
possibility that some of these strings may be infinite.

The recursion for the blocks is
\beql{defB0}
B_1^{(m)} = m \, , 
\eeq
and, for $n \ge 1$,
\beql{defB}
B_{n+1}^{(m)} =  \left\{
	\begin{array}{ll}
	(B_n^{(m)})^{m+1}  S_n^{(m)} & \mbox{~if~} l(B_n^{(m)}) < \infty \, , \\
	B_n^{(m)} & \mbox{~if~} l(B_n^{(m)}) = \infty \, ,
	\end{array}
                 \right.
\eeq
where $S_n^{(m)}$ will be constructed from $B_n^{(m)}$. 
If $l(B_n^{(m)}) = \infty$, $S_i^{(m)} = \emptyset$ for $i \ge n$.
If $l(B_n^{(m)}) < \infty$, consider the sequence
$s_n^{(m)}(1),s_n^{(m)}(2), s_n^{(m)}(3),\ldots \in \PP_m^{\ast}$ defined by
\begin{eqnarray}\label{defsmn}
s_n^{(m)}(1) & = & \sC^{(m)}((B_n^{(m)})^{m+1}) \, , \nonumber \\
s_n^{(m)}(i+1) & = & \sC^{(m)}((B_n^{(m)})^{m+1} s_n^{(m)}(1) \ldots s_n^{(m)}(i))
\quad\mbox{for} \quad i \ge 1 \,.
\end{eqnarray}
Clearly $s_n^{(m)}(1) \ge m+1$.  
If there is an integer $i \ge 1$ such that $s_n^{(m)}(i+1) < m+1$, 
choose the smallest such $i$, and set
\beql{defS}
S_n^{(m)} = s_n^{(m)}(1), s_n^{(m)}(2), \ldots, s_n^{(m)}(i) \in \PP_{m+1}^+ \, ,
\eeq
but if no such $i$ exists set
\beql{defS2}
S_n^{(m)} = s_n^{(m)}(1), s_n^{(m)}(2), \ldots \in \PP_{m+1}^{\ast} \, .
\eeq
In the latter case $S_n^{(m)}$ and $B_{n+1}^{(m)}$ are infinite.

The $T_n^{(m)}$ are defined as follows.
For $n \ge 1$,
if $S_1^{(m)}, \ldots, S_n^{(m)}$ are finite we set
\beql{defT}
T_{n+1}^{(m)} = S_1^{(m)} \ldots S_n^{(m)} \in \PP_{m+1}^{+} \, ,
\eeq
while if $S_1^{(m)}, \ldots, S_{n-1}^{(m)}$ are finite
but $S_n^{(m)}$ is infinite, we still use \eqn{defT} and define
\beql{defT2}
T_{i}^{(m)} = T_{n+1}^{(m)}
\eeq
for $i \ge n+2$.
In the latter case all the $T_{i}^{(m)}$ for $i \ge n+1$ are infinite.
Note that $T_1^{(m)}$ is always undefined.

The lengths of these strings (which may be infinite)
are denoted by
\begin{eqnarray}
\bm(n) & = & l(B_n^{(m)}) \label{defbeta} \, , \\
\sm(n) & = & l(S_n^{(m)}) \label{defsigma} \, , \\
\tm(n) & = & l(T_n^{(m)}) \label{deftau} \,.
\end{eqnarray}
We also let $B^{(m)}= b^{(m)}(1), b^{(m)}(2), b^{(m)}(3), \ldots =
\lim_{n \to \infty} B^{(m)}_n$. This is well defined 
since each $B^{(m)}_n$ starts with $B^{(m)}_{n-1}$.

We will require three lemmas.


\medskip
\begin{lemma}\label{Lem1}
For $m \ge 1$, if $\Am$ contains a string $U^{t+1} \in \PP_m^+$ for some $t \ge m$,
then $U \in \PP_t^+$.
\end{lemma}

\noindent{\bf Proof.} If $t = m$ the claim is trivially true,
so we may assume $t \ge m+1$. 
Suppose, on the contrary, that $U \notin \PP_t^+$.  
Then we may write $U=GiH$ for $G,H \in \PP_m^\ast$ and some $i$
with $m \le i \le t-1$.  
Thus $A$ contains $GiH \, GiH \, \ldots \,
GiH$ $(t+1$ copies).  But the final $i$ is preceded by $t$ copies of
$iHG$, so the final $i$ must be at least $t$, by definition of $\Am$,
a contradiction.~~~$\bsq$


\medskip
\begin{lemma}\label{Lem2}
For $m \ge 1$, $n \ge 2$,
(a) $T_n^{(m)}$ is a suffix of $B_n^{(m)}$, and (b)
this is the only
occurrence of $T_n^{(m)}$ as a substring of $B_n^{(m)}$.
\end{lemma}

\noindent{\bf Proof.} Fix $m \ge 1$.
It follows by iterating \eqn{defB} that
\begin{eqnarray}\label{expandB}
B_{n+1}^{(m)}
& = &
(B_n^{(m)})^m (B_{n-1}^{(m)})^m \cdots (B_1^{(m)})^m B_1^{(m)} 
S_1^{(m)} S_2^{(m)} \ldots  S_{n-1}^{(m)} S_n^{(m)}  \nonumber \\
& = &
(B_n^{(m)})^m (B_{n-1}^{(m)})^m \cdots (B_1^{(m)})^m B_1^{(m)} T_{n+1}^{(m)} \,,
\end{eqnarray}
provided all of $S_1^{(m)} S_2^{(m)} \ldots S_{n}^{(m)}$ are finite.
If $S_1^{(m)} S_2^{(m)} \ldots S_{n-1}^{(m)}$ are finite
but $S_n^{(m)}$ is infinite, \eqn{expandB} is still true, but
\beql{expandB2}
B_{i}^{(m)} =B_{n+1}^{(m)}, ~
T_{i}^{(m)} = T_{n+1}^{(m)}
\mbox{~for~} i \ge n+2 \,.
\eeq
Assertion (a) follows at once.
To show (b) we use induction on $n$.  The base case,
$n=2$, is true because $T_2^{(m)}=m+1$ and $B_2^{(m)} = (m^{m+1},m+1)$.
If $T_{n+1}^{(m)}$ is infinite and has two occurrences
in $B_{n+1}^{(m)}$, they are both suffixes of $B_{n+1}^{(m)}$,
implying that $T_{n+1}^{(m)}$ is a suffix of itself,
and hence is a periodic sequence.
But this is impossible: let $M_0$ be the maximal
element of $T_{n+1}^{(m)}$.  After sufficiently many terms
the curling number given by \eqn{defsmn}
would produce a term greater exceeding $M_0$, a contradiction.
On the other hand, suppose that all the $S_n^{(m)}$ are finite.
If $T_{n+1}^{(m)}$ also occurs in $B_{n+1}^{(m)}$
other than as a suffix, it must be
a substring of a block $B_j^{(m)}$ in \eqn{expandB},
for some $j$ with $2 \le j \le n$, for otherwise it would contain the
$m$ at the beginning of a block.  Write $B_j^{(m)} = U T_{n+1}^{(m)}
V = U T_j^{(m)} S_j^{(m)} \cdots S_n^{(m)} V$ for some $U \in \PP_m^+$ and
$V \in \PP_m^\ast$.  But $l(S_j^{(m)} \cdots S_n^{(m)}) > 0$, so $T_j^{(m)}$
occurs as a non-suffix in $B_j^{(m)}$, a contradiction
to the induction hypothesis.~~~$\bsq$

\noindent{\bf Remark.}
It follows from the above proof that, for any $r$
with $1 \le r \le m+1$,
any finite substring $(B_n^{(m)})^r$ in \eqn{expandB}
contains exactly $r$ copies of $T_n^{(m)}$,
each one occurring at the end of a $B_n^{(m)}$.  The copies are disjoint.


\medskip
\begin{lemma}\label{lemF1}
For $m\ge 1$, $n \geq 2$, suppose that $k=b^{(m)}(i) \geq m+1$ with
$1 \leq i \leq \beta^{(m)}(n)$. Then there exists a $Y$ such that 
$b^{(m)}(1),\ldots, b^{(m)}(i-1) = XY^k$. Moreover, let 
Y satisfy this condition with
l(Y) minimal and suppose $m \in Y$.  Then $Y = B_j^{(m)}$ for some $j$ with $1\leq j
\leq n-1$.
\end{lemma}

\noindent {\bf Proof.}
We fix $m\geq 1$, and will prove the result
for all $n$ by induction.
The base case $n=2$ is immediate,
since $B^{(m)}_2 = (m^{m+1}, m+1)$.
Supposing the result holds for some $n \ge 2$,
we will show it holds for $n+1$.
 If $B^{(m)}_n$ is infinite then
the result also holds for $n+1$, by \eqn{defB}, so we may assume that
$B^{(m)}_n$ is finite. 
Then $B^{(m)}_{n+1} = (B^{(m)}_n)^{(m+1)} S_n^{(m)}$, by \eqn{defB}.
We must show that the result holds for all positions
$\beta^{(m)}(n) <i \leq \beta^{(m)}(n+1)$.

If $i$ is a 
position in $(B_n^{(m)})^{m+1}$, we may write
$i = r\beta^{(m)}(n) + s \leq (m+1)\beta^{(m)}(n)$, for
$1 \le r \le m$, $1\leq s \leq \beta^{(m)}(n)$. 
Then $b^{(m)}(i) = b^{(m)}(s)$ and by induction 
we know that in the first $B_n^{(m)}$ we can write $b^{(m)}(1),\ldots,
b^{(m)}(s-1) = XY^k$, and if the minimal $Y$ contains an $m$ then it
equals $B_j^{(m)}$ for some $j$ with $1 \leq j \leq n-1$. Therefore this $Y$
(and no shorter string) can also be used at position $i$, and thus the statement holds.

If $i = (m+1)\beta^{(m)}(n) + 1$ then the part preceding $i$ is 
$(B_n^{(m)})^{m+1}$, and from \eqn{defsmn} we have
$$
k=b^{(m)}(i) = s_n^{(m)}(1) = \sC^{(m)}((B_n^{(m)})^{m+1}) \geq m+1 \, .
$$
So certainly one $Y$ exists with 
$b^{(m)}(1),\ldots, b^{(m)}(i-1) = XY^k$.
We must show 
that if the minimal $Y$ satisfying this property contains an $m$,
then that $Y=B^{(m)}_j$ for some $j \le n$.
If $Y$ contains an $m$, then it contains $T^{(m)}_{n}$
as a substring, since the last $m$ in $B^{(m)}_{n+1}$ occurs before the
$T^{(m)}_{n}$ in the last copy of $B^{(m)}_n$. 
Therefore the string $Y^k$ contains at least $k \ge m+1$
copies of $T^{(m)}_n$.
It follows from 
the Remark below
Lemma~\ref{Lem2} that $k=m+1$ and $Y=B^{(m)}_n$.

If $i > (m+1)\beta^{(m)}(n) +1$ we see by the definition of $S_n^{(m)}$ that
again a $Y$ exists.
If $Y$ contains an $m$, then it must properly contain 
the $T^{(m)}_n$ in the final copy of $B^{(m)}_n$.
But in the last
of the $k$ copies of
$Y$ the copy of $T^{(m)}_n$ is followed by an
integer larger than $m$, whereas in the earlier $k-1$ copies it was followed by
the first element of $B^{(m)}_n$, which is $m$. 
This is a contradiction, and shows that in this case
$Y$ cannot contain an $m$.~~~$\bsq$

Note that, by definition of $S^{(m)}_{n-1}$, the $Y^k$ for the first element of
$S^{(m)}_n$ goes back further than $B^{(m)}_n$, and thus contains an $m$. 
Therefore we see that the situation described in the penultimate paragraph of the
above proof is indeed the case and we may conclude that
\beql{eqsmn1}
s^{(m)}_{n}(1) = m+1 \mbox{~for~all~} m \ge 1, n \ge 1 \, .
\eeq

At this point we can already see that the concatenation of the glue strings is 
equal to the next $A$ sequence:


\medskip
\begin{theorem}\label{lemF2}
Suppose $m \ge 1$. For all $n\geq 2$,
$T_n^{(m)}$ is a prefix of $A^{(m+1)}$, or is 
all of $A^{(m+1)}$ if $T_n^{(m)}$ is infinite.
\end{theorem}

\noindent {\bf Proof.} 
Again we fix $m$ and use induction on $n$. For $n=2$ the
result is trivial.
Supposing the result holds for some $n \ge 2$, we will show it holds for $n+1$.
If $T_n^{(m)}$ is infinite then clearly the result holds
for $T_{n+1}^{(m)}$, so assume that $T_n^{(m)}$ is finite. 

Write $(B_{n}^{(m)})^{m+1} = U T_{n}^{(m)}$ for some $U \in P_m^+$.
We know that $S_{n}^{(m)}$ begins with $m+1 = s_{n}^{(m)}(1) = \sC^{(m)}(U
T_{n}^{(m)}) = \sC^{(m+1)}(U T_{n}^{(m)}) = \sC^{(m+1)}(T_{n}^{(m)})$.  
The last equality
holds because dropping the $U$ can only decrease the value, but it is
already equal to its minimal value of $m+1$. 
By the induction hypothesis, $T_{n}^{(m)}$ is a prefix of $A^{(m+1)}$,
and therefore $T_{n}^{(m)} s_{n}^{(m)}(1)$ is a prefix of $A^{(m+1)}$. 
For $i \ge 1$, as long as $s_{n}^{(m)}(i) \ge m+1$, we have
\begin{eqnarray*}
s_{n}^{(m)}(i+1)
& = & \sC^{(m)}   (U T_{n}^{(m)} s_{n}^{(m)}(1) \ldots s_{n}^{(m)}(i)) \\
& = & \sC^{(m+1)} (U T_{n}^{(m)} s_{n}^{(m)}(1) \ldots s_{n}^{(m)}(i)) \\
& = & \sC^{(m+1)} (  T_{n}^{(m)} s_{n}^{(m)}(1) \ldots s_{n}^{(m)}(i)) \,.
\end{eqnarray*}
The second equality holds because $s_{n}^{(m)}(i) \ge m+1$.  The third
equality holds because $Y^k$ for $s_{n}^{(m)}(i+1)$ goes back no further
than the beginning of $T_{n}^{(m)}$, as we saw in the 
proof
of the previous lemma.  
Hence $T_{n+1}^{(m)} = T_{n}^{(m)} S_{n}^{(m)}$ is a prefix of $A^{(m+1)}$,
as required.~~~$\bsq$


\begin{theorem}\label{Th1}
For all $m \geq 1$, the sequences $A^{(m)}$ and $B^{(m)}$ coincide.
\end{theorem}

\noindent{\bf Proof.}
Fix $m \ge 1$.
We will show by induction on $n$ that, for all $n \ge 1$, 
$B_n^{(m)}$ is a prefix of $A^{(m)}$
or is all of $A^{(m)}$ if $B_n^{(m)}$ is infinite.
This will establish the theorem.

\medskip
The cases $n=1$ and $n=2$ are immediate, 
since $B_1^{(m)}=m$, $B_2^{(m)} = (m^{m+1},m+1)$.
So assume the truth of the induction
hypothesis up to and including some $n \ge 2$.


If $B_n^{(m)}$ is infinite the result follows from \eqn{defB},
so we may assume that $B_n^{(m)}$ and hence $T_n^{(m)}$ are finite.
We wish to show that
$B_{n+1}^{(m)} = (B_n^{(m)})^{m+1} S_n^{(m)}$ is a prefix of $\Am$.
If this is not true, the 
first discrepancy between $B_{n+1}^{(m)}$ and $\Am$
occurs in the substring $(B_n^{(m)})^{m+1}$, by the definition of $S_n^{(m)}$. 
Let $i \ge \bm(n)+1$ be the first
position in $(B_n^{(m)})^{m+1}$ at which $a^{(m)}(i) \ne b^{(m)}(i)$.
Our goal is to show that the existence of $i$ leads to a contradiction.

We may write $i = j \bm(n)+r$ with $1 \le j \le m$ and $1 \le r
\le \bm(n)$. Then $i$ is also minimal with respect to the condition that
$a^{(m)}(i) \ne a^{(m)}(r)$.  
Let $a^{(m)}(1), \ldots, a^{(m)}(i-1)
= X Y^k$ with $k$ maximal and $l(Y)$ minimal.  Then $a^{(m)}(i) = $
max$\{m,k\}$.

We consider two cases, depending on whether or not $a^{(m)}(i)$ is at
the beginning of one of the $\Bnm$ blocks, i.e. whether $r=1$ or $r \ge 2$.  

First, suppose $r=1$;
then we need to prove that $a^{(m)}(i) = a^{(m)}(1) = m$.  
This follows by definition of $S_{n-1}^{(m)}$ if $j=1$,
so assume $j \ge 2$, and that $k = a^{(m)}(i) \ge m+1$.
Using \eqn{expandB} we may write $a^{(m)}(1), \dots, a^{(m)}(i-1) =
(B_n^{(m)})^j = (B_n^{(m)})^{j-1} U m T_n^{(m)}$ for some $U \in \PP_m^\ast$.  
If $T_n^{(m)}$ is a proper suffix of $Y^k$ then $m \in Y$,
which implies that $T_n^{(m)}$ is
a proper suffix of $Y$ and therefore $(B_n^{(m)})^j$ contains at least $m+1$
copies of $T_n^{(m)}$, contradicting the Remark following Lemma~\ref{Lem2}.
On the other hand, if $Y^k$ were a suffix of $T_n^{(m)}$, this would contradict
the fact that $S_{n-1}^{(m)}$ is followed by an element $\le m$. 

Second, suppose that $r \ge 2$.
Let $L = a^{(m)}(1), \ldots, a^{(m)}(r-1)$ and write $L = X_\ast Y_\ast^{k_\ast}$
with $k_\ast$ maximal and $l(Y_\ast)$ minimal.  Then $a^{(m)}(r) = $max$\{m,k_\ast\}$.
By the definition of $i$, $a^{(m)}(i) > a^{(m)}(r) \ge m$. Hence $a^{(m)}(i)
= k \ge m+1$.  To have $a^{(m)}(i) > a^{(m)}(r)$, $L$ must be a suffix
of $Y^k$, so $m \in Y$ and therefore, by Lemma~\ref{Lem1}, $k$ is at most $m+1$
and therefore is equal to $m+1$.  Hence $k_\ast \le m$.

The situation, then, is that $(B_n^{(m)})^j L$ is a prefix of $\Am$.  
We are supposing that we can achieve $a^{(m)}(i) = m+1$ by allowing
$L$ to be a suffix of $Y^{m+1}$.  Noting that $T_n^{(m)}$ is a suffix of
$(B_n^{(m)})^j$, by \eqn{expandB}, we distinguish two cases,
depending on the relationship between $T_n^{(m)} L$ and $Y^{m+1}$.

(i) Suppose that $T_n^{(m)} L$ is a suffix of $Y^{m+1}$.  We know $m \in Y$ and
$m \not\in T_n^{(m)}$, so $Y^{m+1}$ contains at least $m$ disjoint copies
of $T_n^{(m)}$.
Hence $j=m$, and there are exactly $m$ disjoint copies, by
the Remark following Lemma~\ref{Lem2}.
This means that each copy of $T_n^{(m)}$ straddles the end of
one copy of $Y$ and the beginning of the next 
(if not, $T_n^{(m)}$ is wholly contained in $Y$,
and so there are $m+1$ copies of $T_n^{(m)}$ in the
sequence before position $i$, which is
a contradiction since there are only $m$ copies, one in each
of the $m$ copies of $B_n^{(m)}$ and none so far in the next
copy of $B_n^{(m)}$ that we are building),
and hence that $Y$ is a
proper suffix of $T_n^{(m)} L$.  
Write $T_n^{(m)} = VW$ where
$W$ is the intersection of $T_n^{(m)}$ and the last 
(or $(m+1)$-st) copy of $Y$,
and write $B_n^{(m)} = UT_n^{(m)}$, using \eqn{expandB}.  
If $m \ge 2$ it is easy to complete the proof.
We have $Y = WL = WUV$, so $L = UV$
and therefore $i > l(Y^{m+1}) = (m+1)l(WUV) =
(m+1)l(UVW) = (m+1)\bm(n)$, contradicting the definition of $i$.

Suppose then that $m=1$.
Again $L$ is a proper suffix of $Y$ and $Y$ is a proper 
suffix of $T_n^{(1)} L$. Write $Y=WL$, and let $s\geq 2$ be the first 
element of $Y$.
Let this element $s$ in the second copy of $Y$ be preceded by $s$ 
copies of some string $Y'$ with $l(Y')$ minimal.

Suppose that $Y'$ does not contain a $1$.
Since $Y$ does contain a $1$ ($L$ starts with a $1$),
$Y'^s$ is a suffix of the 
first copy of $Y$, and hence also of the second copy of $Y$.
This contradicts the minimality of $l(Y)$,
since then $l(Y')<l(Y)$.

So we may assume that $1\in Y'$, hence by Lemma~\ref{lemF1} we know that
$$
Y'=B_{\kappa}^{(1)}
$$
for some ${\kappa}<n$.
$T_{\kappa}^{(1)}$ is a suffix of $Y^2$ and since $L$ starts with a $1$,
$T_{\kappa}^{(1)}$ is also a suffix of $L$.
By Lemma~\ref{Lem2}, $B_{\kappa}^{(1)}$ is also a suffix of $L$.
Suppose $B_{\kappa}^{(1)}=L$.
Then $WL$ is a suffix of $Y'Y'=LL$
(look at the first copy of $Y = WL$ and remember $Y'Y'$ begins with a $1$)
and hence $W$ is a suffix of $L$. But 
then $W^2$ is a suffix of $T_n^{(1)}$, contradicting the fact that $L$ 
starts with a $1$.

So we may assume that $B_{\kappa}^{(1)}$ is a strict suffix of $L$.
But now $l(Y)> l(L)\geq 2 l(B_{\kappa}^{(1)})$.
(Indeed, if $l(L)<2l(B_{\kappa}^{(1)})$,
then we know that $L$ is a prefix of $B_n^{(1)}$, by definition,
$B_{\kappa}^{(1)}B_{\kappa}^{(1)}$ is also a prefix of $B_n^{(1)}$,
and so $L$ is a strict prefix of $B_{\kappa}^{(1)}B_{\kappa}^{(1)}$;
but $L$ has $T_{\kappa}^{(1)}$ as a suffix, so by Lemma~\ref{Lem2}, 
$L=B_{\kappa}^{(1)}$, a contradiction.)
But now $(B_{\kappa}^{(1)})^2$ is a suffix of $Y$, contradicting the minimality 
of $Y$.

(ii) Suppose on the other hand
that $Y^{m+1}$ is a suffix of $T_n^{(m)} L$.  Since no $Y$ is
contained in $T_n^{(m)}$ (remember that $m \in Y$),
$Y^m$ is a suffix of $L$ and the
first element, $t$, of $Y$ is in $T_n^{(m)}$ with $t \ge m+1$.  Therefore the
first element of the second $Y$ is also $t$ and since $(B_n^{(m)})^j L$
is a prefix of $\Am$, $Y$ ends with $U^t$ for some $U$.  Hence $U^t$
is a suffix of $L$, which contradicts the fact that $k_\ast = m$.
This completes the proof.~~~$\bsq$


\begin{coro}\label{Cor2}
The sequence $A^{(m)}$ contains every integer $\ge m$.
\end{coro}

\noindent {\bf Proof.}
 From Theorem \ref{lemF2} we know that, for $m\ge 2$, $n \ge 2$,
$T_{n}^{(m-1)}$ is a prefix of $A^{(m)}$, so, for a given $m$, either
$$
A^{(m)} = S_1^{(m-1)}  S_2^{(m-1)} \ldots S_{n}^{(m-1)}
$$
if some  $S_n^{(m-1)}$ is infinite, or
$$
A^{(m)} = S_1^{(m-1)}  S_2^{(m-1)} S_3^{(m-1)} \dots
$$
if all  $S_n^{(m-1)}$ are finite.
Also, by Theorem \ref{Th1}, $B_{n+1}^{(m-1)}$ is a prefix of $A^{(m-1)}$,
so from \eqn{defB},
if some  $S_n^{(m-1)}$ is infinite, $A^{(m-1)}$ contains 
$$
S_1^{(m-1)},  S_2^{(m-1)}, \ldots, S_{n}^{(m-1)} \,,
$$
or if all  $S_n^{(m-1)}$ are finite, $A^{(m-1)}$ contains
$S_1^{(m-1)},  S_2^{(m-1)}, \ldots, S_{n}^{(m-1)}$
for all $n$.
In either case (and this is the key point), every prefix of
$A^{(m)}$ is a subsequence of $A^{(m-1)}$.
Repeating this argument shows that every prefix of
every $A^{(j)}$ is a subsequence of $A^{(m)}$ if $j\geq m$.

Since $A^{(j)}$ begins with $j$, $A^{(m)}$ contains every 
integer $j \ge m$.~~~$\bsq$


\begin{coro}\label{CorFinite}
The strings $B_n^{(m)}$, $S_n^{(m)}$ and $T_n^{(m)}$ have finite length.
\end{coro}

\noindent {\bf Proof.} 
The first occurrence of an integer in $A^{(m)}$ is necessarily followed 
by an $m$. Since we saw in the previous corollary that $A^{(m)}$ contains
infinitely many different integers, it follows that all $S_n^{(m)}$ are finite.
This implies that $B_n^{(m)}$ and $T_n^{(m)}$ are also finite.~~~$\bsq$

\section*{Promotion}\label{SPROM}

In the definition of $A^{(m)}$, \eqn{defA}, let us say
that $a^{(m)}(i)$ is {\em promoted}
if either $i=1$ or \\
$\sC (a^{(m)}(1),\ldots, a^{(m)}(i-1)) <m$.
If we know which elements in $A^{(m+1)}$ are promoted,
we can recover $A^{(m)}$ from $A^{(m+1)}$.  To make this precise,
we define the 
strings
$D^{(m)}_i \in \PP_m^+$ by
$D^{(m)}_{0} = m$ and, for $i > 0$,
\beql{EqProm1}
D^{(m)}_i = \left\{ \begin{array}{ll}
  D^{(m)}_{i-1} a^{(m+1)}(i) & \mbox{ if $a^{(m+1)}(i)$ is not promoted} \, , \\
 \left( D^{(m)}_{i-1} \right)^{m+1} a^{(m+1)}(i) & \mbox{ if $a^{(m+1)}(i)$
 is promoted} \, .  \end{array} \right.
\eeq
Since $D^{(m)}_i$ starts with $D^{(m)}_{i-1}$, we can define the limiting
sequence $D^{(m)} = \lim_{i \to \infty} D^{(m)}_i$.
Then it can be shown that:
\begin{theorem}
For all $m \ge 1$, the sequences $A^{(m)}$ and $D^{(m)}$ coincide.
\end{theorem}

We omit the proof, which involves arguments similar to those 
used to prove Theorems \ref{lemF2} and \ref{Th1}.
The main difference is that this proof does not require
the finiteness of the glue strings $S_n^{(m)}$. 
Furthermore,
the glue strings now by definition unite to form the next $A$-sequence, but 
on the other hand it becomes more difficult to show that
they are indeed substrings of $A^{(m)}$ itself.


\section{Estimates for the rate of growth}\label{S4}
In this section we take an experimental approach,
and record a series of observations about the sequence.
These observations appear to be correct, 
but we have been unable to prove them.
In \S\ref{SecRuler} we study the lengths of the glue strings $S_n^{(m)}$.
Although these lengths are somewhat irregular, it appears
that they can be ``smoothed'' so as to become 
much more regular ``ruler'' sequences,
whose peak values will be denoted by $\rho^{(m)}(n)$.
In \S\ref{Sec4BB} we describe a ``tabular'' construction for
the higher-order sequences $A^{(2)}, A^{(3)}, \ldots$
which leads to a recurrence relating
the $\rho^{(m)}(n)$, $\beta^{(m)}(n)$ and $\sigma^{(m)}(n)$.
Sections \ref{Sec4B}, \ref{Sec41} and \ref{Sec42}
contain estimates for $\beta^{(m)}(n)$, $\rho^{(m)}(n)$ and
$\tau^{(m)}(n)$.
Finally, in \S\ref{Sec43}, we use
these estimates to determine where each number $t \ge 1$ appears for the
first time in our sequence $A$. 

\subsection{Ruler sequences and smoothing }\label{SecRuler}
It appears that the sequence $\sm = \sm(1), \sm(2), \sm(3), \ldots$ 
giving the lengths of the glue strings $S_n^{(m)}$
is essentially
a ``ruler'' sequence, in the sense that 
$\sm(n)$ essentially depends only on
the $(m+1)$-adic valuation of $n$.

For positive integers $m,n$, define the {\em $m$-adic valuation} of $n$,
$|n|_m$, to be the highest power of $m$ dividing $n$.  The classical
example of a ruler sequence is the sequence $r= r(1), r(2), r(3), \ldots$
given by
\beql{EqR9}
r(n) = |n|_2 +1 \,.
\eeq
The first 32 terms are
$$
\begin{array}{cccccccccccccccc}
\bone & \btwo & 1 & \bthree & 1 & 2 & 1 & \bfour & 1 & 2 & 1 & 3 & 1 & 2 & 1 & \bfive \\
1 & 2 & 1 & 3 & 1 & 2 & 1 & 4 & 1 & 2 & 1 & 3 & 1 & 2 & 1 & \bsix
\end{array}
$$
where the new record entries, shown in bold-face, occur at powers of 2.
For much more about this sequence, including an extensive bibliography,
see entry A1511 in \cite{OEIS}.

The initial values of $\sigma^{(1)}, \ldots, \sigma^{(4)}$ are shown
in Table~\ref{T4}, and the record entries in $\sigma^{(1)}, \ldots,
\sigma^{(10)}$ in Table~\ref{T5}.  Let $\pi^{(m)}(j)$ ($j \ge 0 $)
denote the $j$-th record in $\sigma^{(m)}$.

\begin{table}[htb]
$$
\begin{array}{ccccccccccccccccc}
n & 1 & 2 & 3 & 4 & 5 & 6 & 7 & 8 & 9 & 10 & 11 & 12 & 13 & 14 & 15 & 16 \\
\sigma^{(1)} (n) & \bone & \bthree & 1 & \bnine & 4 & \btf & 1 & 3 & 1 & 9 & 4 & \bss & 1 & 3 & 1 & 9 \\
\sigma^{(2)} (n) & \bone & 1 & \bthree & 1 & 1 & 3 & 1 & 1 & \bnine & 1 & 1 & 3 & 1 & 1 & 3 & 1 \\
\sigma^{(3)} (n) & \bone & 1 & 1 & \bthree & 1 & 1 & 1 & 3 & 1 & 1 & 1 & 3 & 1 & 1 & 1 & \bten \\
\sigma^{(4)} (n) & \bone & 1 & 1 & 1 & \bthree & 1 & 1 & 1 & 1 & 3 & 1 & 1 & 1 & 1 & 3 & 1 \\ [+.25in]
n & 17 & 18 & 19 & 20 & 21 & 22 & 23 & 24 & 25 & 26 & 27 & 28 & 29 & 30 & 31 & 32 \\
\sigma^{(1)} (n) & 4 & 24 & 1 & 3 & 1 & 9 & 4 & \bons & 3 & 1 & 9 & 4 & 24 & 1 & 3 & 1 \\
\sigma^{(2)} (n) & 1 & 9 & 1 & 1 & 3 & 1 & 1 & 3 & 1 & 1 & \btt & 1 & 3 & 1 & 1 & 3 \\
\sigma^{(3)} (n) & 1 & 1 & 1 & 3 & 1 & 1 & 1 & 3 & 1 & 1 & 1 & 3 & 1 & 1 & 1 & 10 \\
\sigma^{(4)} (n) & 1 & 1 & 1 & \beleven & 1 & 1 & 1 & 1 & 3 & 1 & 1 & 1 & 1 & 3 & 1 & 1  \\ [+.25in]
n & 33 & 34 & 35 & 36 & 37 & 38 & 39 & 40 & 41 & 42 & 43 & 44 & 45 & 46 & 47 & 48 \\
\sigma^{(1)} (n) & 9 & 4 & 68 & 3 & 1 & 9 & 4 & 24 & 1 & 3 & 1 & 9 & 4 & \bfeo & 3 & 1
\end{array}
$$

\caption{Values of $\sigma^{(1)}(n)$ for $n \le 48$ and 
$\sigma^{(2)} (n)$, $\sigma^{(3)}
(n)$, $\sigma^{(4)} (n)$, for $n \le 32$, with record entries shown
in bold-face.}
\label{T4}
\end{table}

\medskip

\begin{table}[htb]
$$
\begin{array}{c|rrrrrrrrrr}
m \setminus j & \multicolumn{1}{r}{0} &
\multicolumn{1}{r}{1} &
\multicolumn{1}{r}{2} &
\multicolumn{1}{r}{3} &
\multicolumn{1}{r}{4} &
\multicolumn{1}{r}{5} &
\multicolumn{1}{r}{6} &
\multicolumn{1}{r}{7} &
\multicolumn{1}{r}{8} &
\multicolumn{1}{r}{9} \\ \hline
1 & 1 & 3 & \mathit{9} & 24 & 67 & \mathit{196} & \mathit{581} & \mathit{1731} & \mathit{5180} & \mathit{15534} \\
2 & 1 & 3 & 9 & \mathit{32} & 119 & 463 & 1837 & \mathit{7332} & 29307 & 117203 \\
3 & 1 & 3 & 10 & 42 & \mathit{200} & 983 & 4892 & 24434 & 122141 \\
4 & 1 & 3 & 11 & 55 & 315 & \mathit{1872} & 11205 & 67195 \\
5 & 1 & 3 & 12 & 70 & 471 & 3273 & \mathit{22883} \\
6 & 1 & 3 & 13 & 87 & 673 & 5355 & 42805 \\
7 & 1 & 3 & 14 & 106 & 927 & 8309 & 74740 \\
8 & 1 & 3 & 15 & 127 & 1239 & 12351 & 123463 \\
9 & 1 & 3 & 16 & 150 & 1615 & 17721 \\
10 & 1 & 3 & 17 & 175 & 2061 & 24683
\end{array}
$$

\caption{
Values of $\pi^{(m)} (j)$, the $j$-th record in sequence
$\sigma^{(m)}$.  
The smoothed record values $\rho^{(m)} (j)$ are obtained by reducing
the italicized entries by $1$.
The next three terms in the first row are $\mathit{46578}$,
$\mathit{139713}$, $\mathit{419116}$, and the next term in the $m=2$ row is 468785.
The missing entries in this table have not been calculated, although we predict that
the entries on or below the diagonal $m=j$ are given by \eqn{EqR16}
and the entries just above this diagonal by \eqn{EqR16a}.
}
\label{T5}
\end{table}

\newpage

As can be seen from Table~\ref{T4}, $\sigma^{(1)}$ is not
quite as regular as the ruler sequence $r$.  However:
\begin{obs}\label{Obs1}
If the sequence $\sigma^{(1)}$ is ``smoothed''
by replacing every instance of $4$ by the pair of numbers $3$, $1$, every
$9$ by $8$, $1$, every $25$ by $24$, $1$, and so on, $\sigma ^{(1)}$ becomes 
a ruler sequence $r^{(1)}$ given by
\beql{EqR9a}
r^{(1)} (n) = \rho^{(1)} (|n|_2) \,,
\eeq
in which the first $64$ terms are
$$
\begin{array}{cccccccccccccccc}
\bone & \bthree & 1 & \beight & 1 & 3 & 1 & \btf & 1 & 3 & 1 & 8 & 1 & 3 & 1 & \bss \\
1 & 3 & 1 & 8 & 1 & 3 & 1 & 24 & 1 & 3 & 1 & 8 & 1 & 3 & 1 & \bonf \\
1 & 3 & 1 & 8 & 1 & 3 & 1 & 24 & 1 & 3 & 1 & 8 & 1 & 3 & 1 & 67 \\
1 & 3 & 1 & 8 & 1 & 3 & 1 & 24 & 1 & 3 & 1 & 8 & 1 & 3 & 1 & \bfeh 
\end{array}
$$
and where the record values (shown in bold-face)
$\rho^{(1)} (0), \rho^{(1)} (1), \ldots$ are
\beql{EqR10}
1,3,8,24,67,195,580, 1730, 5179, 15533,
46578, 139712, 419115, \ldots\,.
\eeq
The numbers $i$ in $\sigma^{(1)}$ that are to be replaced by $i-1$, $1$
to get $r^{(1)}$ are
\beql{EqR11}
4, 9, 25, 68, 196, 581, 1731, 5180, 15534, 46579, 139713, 419116, \ldots \,,
\eeq
\end{obs}

The numbers that need to be smoothed, given
in \eqn{EqR11}, are one greater than the numbers in \eqn{EqR10},
except that $2$ is missing.
The records in the smoothed sequence $r^{(1)}$, \eqn{EqR10},
either agree with or are one
less than the terms in the first row of Table~\ref{T5}.

The sequences $\sigma^{(m)}$ for $m \ge 2$ appear to need less smoothing
than $\sigma^{(1)}$ to make them into ruler sequences.  In the range of
our tables, $\sigma^{(2)}$ needs to be smoothed by replacing every 32
by 31, 1, and every 7332 by 7331, 1; $\sigma^{(3)}$ by replacing every
200 by 199, 1; $\sigma^{(4)}$ by replacing every 1872 by 1871, 1; 
and so on.
If $r^{(m)}$ denotes the smoothed version of $\sigma^{(m)}$
and $\rho^{(m)}(j)$ the $j$-th record in the
smoothed version (see Table~\ref{T5}) then we have, for all $m \ge 1$, $n \ge 1$,
\beql{EqR12}
r^{(m)} (n) = \rho^{(m)} (|n|_{m+1})\,.
\eeq

The lengths $\beta^{(m)} (n)$ of the blocks are given by
(from \eqn{defB}, \eqn{defbeta})
\begin{eqnarray}\label{EqR17a}
\beta^{(m)}(1)  & = & 1 \, , \nonumber \\
\bm (n+1)  & = &  (m+1) \bm (n) + \sm (n) \mbox{~for~} n \ge 1 \,.
\end{eqnarray}

The initial values of 
$\beta^{(1)}(n), \ldots, \beta^{(6)}(n)$
are shown in Table~\ref{TBlen}.

\begin{table}[htb]
$$
\begin{array}{c|rrrrrrrrr}
m \setminus n & \multicolumn{1}{r}{1} &
\multicolumn{1}{r}{2} &
\multicolumn{1}{r}{3} &
\multicolumn{1}{r}{4} &
\multicolumn{1}{r}{5} &
\multicolumn{1}{r}{6} &
\multicolumn{1}{r}{7} &
\multicolumn{1}{r}{8} \\ \hline
1 & 1 & 3 &  9 & 19 & 47 & 98 & 220 & 441 \\
2 & 1 & 4 & 13 & 42 & 127 & 382 & 1149 & 3448 \\
3 & 1 & 5 & 21 & 85 & 343 & 1373 & 5493 & 21973 \\
4 & 1 & 6 & 31 & 156 & 781 & 3908 & 19541 & 97706 \\
5 & 1 & 7 & 43 & 259 & 1555 & 9331 & 55989 & 335935 \\
6 & 1 & 8 & 57 & 400 & 2801 & 19608 & 137257 & 960802 
\end{array}
$$
\caption{
Lengths $\beta^{(m)}(n)$ of the blocks $B^{(m)}_n$.}
\label{TBlen}
\end{table}

\subsection{The tabular construction }\label{Sec4BB}

The appearance of ruler sequences can be partially explained 
if we present the construction of the higher-order
sequences $A^{(2)}, A^{(3)}, \ldots$,
in a tabular format.
In this construction we keep track not only of the actual
value $A^{(m)}(n) = \max \{ m, k\}$ (cf. \eqn{defA})
but also whether the promotion rule was invoked (if 
$k < m$ we indicate this by drawing a circle around
the entry) and the length of the shortest $Y$
that was used to compute $k$ if $k \ge m$ (shown
as a subscript; if the promotion
rule was invoked the subscript is $0$).
This tabular construction will also suggest a recurrence 
that relates 
$\rho^{(m)}(n+1)$, $\beta^{(m+1)}(n+1)$ and $\sigma^{(m+1)}(n+1)$.

We will construct $A^{(2)}$ as an example.
We start by making a small table of the glue strings $S^{(2)}_n$
for $n \le 10$ --- see Table~\ref{TabBB1}.
(We already saw $S_1^{(2)}, \ldots, S_4^{(2)}$ in Table~\ref{T3}.)

\begin{table}[htb]
$$
\begin{array}{c|l}
n & S^{(2)}_n \\ \hline
1 & 3 \\
2 & 3 \\
3 & 3 ~ 3 ~ 4 \\
4 & 3 \\
5 & 3 \\
6 & 3 ~ 3 ~ 4 \\
7 & 3 \\
8 & 3 \\
9 & 3 ~ 3 ~ 4 ~ 3 ~ 3 ~ 3 ~ 3 ~ 4 ~ 4 \\
10 & 3 
\end{array}
$$
\caption{ The first few glue strings $S^{(2)}_n$.}
\label{TabBB1}
\end{table}

We know from Section \ref{S3} that
$A^{(2)} =
\displaystyle\lim_{n \to \infty} B_n^{(2)}
= \displaystyle\lim_{n \to \infty} T_n^{(1)}
= S_{1}^{(1)} S_{2}^{(1)} \ldots
$
and that $B_{n+1}^{(2)} = (B_n^{(2)})^3 S_n^{(2)}$.
Table~\ref{TabBB2} shows the beginning of the construction of $A^{(2)}$.

The aim is to understand how $A^{(2)}$ breaks into
the consecutive $S_{n}^{(1)}$ glue strings for $A^{(1)}$.
To do this, a version of $A^{(2)}$ is produced
in which terms that are obtained by promotion are circled,
and where the subscript on each term is either $0$ for
a circled term or else gives the length of the shortest $Y$
that can be used to compute that term.
The circled terms will be the first terms of each of
the glue strings $S_{n}^{(1)}$ of $A^{(1)}$.
Most of the circling and subscripting work is done by a few simple rules.
However, the rules occasionally give the wrong answer
and a few corrections may need to be made by hand
at the end of each round. It is the presence of
these adjustments that makes our sequence hard to analyze.

We start with $B_1^{(2)} = \cto$.
The rules for going from  $B_n^{(2)}$ to $B_{n+1}^{(2)}$ 
are as follows: \\
(i) Write $B_n^{(2)}$ as a single string,
and construct a three-rowed array in which
each row is a copy of $B_n^{(2)}$, omitting
all circles from the third row.
This three-rowed array (after $S^{(2)}_n$ is appended
in step (iii)) will form $B_{n+1}^{(2)}$
when read as a single string.
(When constructing $A^{(m)}$ we make $m$ copies of $B_n^{(m)}$
and omit the circles from the $m$-th copy.) \\
(ii) The subscripts in rows 2 and 3 are the same
as in row 1, except that terms in row 3 that are
under circled terms in row 2 have their
subscripts changed to $l(B_n^{(2)})$. \\
(iii) Append $S^{(2)}_n$ to the end of row 3.
The first term of $S^{(2)}_n$ receives the subscript $l(B_n^{(2)})$.
The subscripts on the remaining terms of $S^{(2)}_n$
must be computed separately---they can be obtained
from the tabular construction of $A^{(m+1)}$. \\
(iv) Finally, a few circles in row 2 may need to
be omitted and their subscripts recomputed,
as well as the subscripts on the same terms
in row 3.

In Table~\ref{TabBB2}, rules (i)--(iii) give the correct
answers for $B_2^{(2)}$ and $B_3^{(2)}$. But in $B_4^{(2)}$ 
four terms (marked with asterisks in Table~\ref{TabBB2}) must be corrected.
The first entry in row 2 of $B_4^{(2)}$ is $\cto$. 
However, row 1 ends with 3 3 = $S_1^{(2)} S_2^{(2)} = Y^2$,
with a $Y$ of length 1,
so that 2 did not need to be promoted and we must
change $\cto$ to $2_1$.
The fifth entry in row 2 of $B_4^{(2)}$ is $\cto$.
But it is preceded by 
$$
3 \, 2 \, 2 \, 2 \, 3 ~~ 3 \, 2 \, 2 \, 2 \, 3  ~=~ S_1^{(2)} B_2^{(2)} S_2^{(2)} B_2^{(2)} ~=~ Y^2 \, ,
$$
with a $Y$ of length 5,
so we must change $\cto$ to $2_5$.
The corresponding entries in row 3, presently both equal to $2_{13}$,
also get changed to $2_1$ and $2_5$ respectively.

\begin{table}[htb]
$$
\begin{array}{lcccccccccccccccc}
B_1^{(2)} = & \cto \\  [+.10in]

B_2^{(2)} = & \cto \\
            & \cto \\
            & 2_1  & 3_1 \\ [+.10in]

B_3^{(2)} = & \cto & \cto & 2_1 & 3_1 \\
            & \cto & \cto & 2_1 & 3_1 \\
            & 2_4  & 2_4  & 2_1 & 3_1 & 3_4 \\ [+.10in]

B_4^{(2)} = & \cto & \cto & 2_1 & 3_1 & \cto & \cto & 2_1 & 3_1 & 2_4 & 2_4 & 2_1 & 3_1 & 3_4 \\
            & \cto^{\ast} & \cto & 2_1 & 3_1 & \cto^{\ast} & \cto & 2_1 & 3_1 & 2_4 & 2_4 & 2_1 & 3_1 & 3_4 \\
            & 2_{13}^{\ast} & 2_{13} & 2_1 & 3_1 & 2_{13}^{\ast} & 2_{13} & 2_1 & 3_1 & 2_4 & 2_4 & 2_1 & 3_1 & 3_4 & 3_{13} & 3_1 & 4_1
\end{array}
$$

\caption{ Tabular construction of $A^{(2)}$. }
\label{TabBB2}
\end{table}

When we extend Table~\ref{TabBB2} to $B_{10}^{(2)}$, we find
that in all only ten circles need to be removed.
After $B_4^{(2)}$, the next 
changes are at $B_7^{(2)}$, where two circles get removed
because of the splittings $S_4^{(2)} S_5^{(2)}$ = 3 3 = $Y^2$,
with a $Y$ of length 1, and
$S_4^{(2)} B_5^{(2)} S_5^{(2)} B_5^{(2)} = Y^2$ with a $Y$ of length 128.
But not all instances of such splittings
cause circles in the table to be removed, and not
all circle-removals arise in this way.
It seems difficult to explain exactly where corrections
to the table are required.

However, the corrections are rare, and still fewer
corrections are needed for larger values of $m$.

Since $A^{(2)}$ is also  
$\displaystyle\lim_{n \to \infty} T_n^{(1)}$,
we can read off the lengths of the glue strings $S_n^{(1)}$ 
from the table. Look at the lengths of the strings
(in $B_4^{(2)}$) between one circle and the next: these are
1, 3, 1, 9, 4, 24, $\ldots$. exactly the values of $\sigma^{(1)}(1)$,
$\sigma^{(1)}(2), \ldots$ (cf. Table~\ref{T4}).
If we do not make the corrections needed in step (iv),
we instead get the smoothed lengths 1, 3, 1, 8, 1, 3, 1, 24, $\ldots$.
These observations lead to our conjectured recurrence.
For example, note that the string in $B_4^{(2)}$
from the last circled entry to the end has length 24
(which is $\rho^{(1)}(3)$)
and is made up of the last string in $B_3^{(2)}$ (length 8,
which is  $\rho^{(1)}(2)$)
plus the whole of $B_3^{(2)}$ (length 13, which is $\beta^{(2)}(3)$),
plus $S_3^{(2)}$ (length 3, which is $\sigma^{(2)}(3)$).
More generally, we have:

\begin{obs}\label{Obs2}
For $m \ge 1$,
\begin{eqnarray}\label{REC1}
\rho^{(m)}(0)  & = & 1 \, , \nonumber \\
\rho^{(m)}(n+1)  & = &  \rho^{(m)}(n) + \beta^{(m+1)}(n+1)
+ \sigma^{(m+1)}(n+1) \mbox{~for~} n \ge 0 \,.
\end{eqnarray}
\end{obs}

This recurrence is supported by all the data, although we
do not have a proof.

\subsection{ Estimates for the lengths $\beta^{(m)} (n)$ of the blocks}\label{Sec4B}
In this section we first prove formulas \eqn{EqBz}, \eqn{EqBa}, \eqn{EqBb},
which give the exact vlaue of $\beta^{(m)}(n)$ in the
parabolic region $1 \le n \le (m+1)^2-1$ for $m \ge 1$.
We then give {\em conjectural} estimates for $\beta^{(m)}(n)$ 
(indicated by $\approx$), \eqn{EqBc} and \eqn{EqBd},
which apply for all $m$ and $n$.

We take \eqn{EqR17a} as our starting point.
For $m \ge 1$, $\sigma^{(m)}(n)$ and $r^{(m)}(n)$
coincide for $1 \le n \le (m+1)^2-1$, and
in that range are given by
\beql{EqBy}
\sigma^{(m)}(n) = r^{(m)}(n) = \left\{ 
\begin{array}{ll}
  1 & \mbox{~if~} m+1 \mbox{~does~not~divide~} n \, , \\
  3 & \mbox{~if~} m+1 \mbox{~divides~} n \, .
\end{array} \right.
\eeq
By iterating \eqn{EqR17a} and using \eqn{EqBy} we find:
 
\begin{lemma}\label{lemB1}
For $m \ge 1$ and $1 \le n \le (m+1)^2-1$,
\beql{EqBz}
\beta^{(m)}(n) ~=~
\frac{(m+1)^n-1}{m} ~+~
2 ~ \frac{(m+1)^{n-1}-(m+1)^v}{(m+1)^{m+1}-1} \, ,
\eeq
where $v \in \{0, 1, \ldots,m\}$ is given by $n-1 \equiv v \bmod m+1$.
\end{lemma}

In particular, for $1 \le n \le m+1$ we have $v = n-1$ and so
\beql{EqBa}
\beta^{(m)}(n) = \frac{(m+1)^n-1}{m} \, ,
\eeq
and for $m+2 \le n \le 2m+2$ we have $v = n-m-2$ and 
\beql{EqBb}
\beta^{(m)}(n) = 
\frac{(m+1)^n+2(m+1)^{n-m-1} -2(m+1)^{n-m-2} -1}{m}  \, .
\eeq
Equation \eqn{EqBa} explains the entries 
on or below the diagonal $n = m+1$ in Table~\ref{TBlen},
\eqn{EqBb} explains the entries in the ``wedge'' $m+2 \le n \le 2m+2$,
and
\eqn{EqBz} the entries in the whole 
parabolic region bounded by $n \le (m+1)^2-1$.

The right-hand side of \eqn{EqBz} is also a good approximation 
to $\beta^{(m)}(n)$ for fixed $m \ge 2$ and $n \ge (m+1)^2$.
The case $m=1$ is special, because of the
greater differences between
$\sigma^{(m)}(n)$ and $r^{(m)}(n)$ when $m=1$.
However, $\beta^{(1)}(n)$ is well-approximated by
\beql{EqBc}
\beta^{(1)}(n) \, \approx \, \varepsilon _1 \, 2^{n-1} \, ,
\eeq
where $\varepsilon _1 = 3.48669886\ldots$.

For our applications, the approximation
\beql{EqBd}
\beta^{(m)}(n) \, \approx \, \varepsilon_m \, (m+1)^{n-1} 
 \mbox{~for~} m \ge 1, n \ge 1
\eeq
(consistent with \eqn{EqBz}--\eqn{EqBc})
will be adequate, where
$\varepsilon_m$ is a constant on the order of 1.

\subsection{ Estimates for the records $\rho^{(m)}(n)$ } \label{Sec41}

We now apply \eqn{REC1} to estimate
$\rho^{(m)} (n)$.
Except for \eqn{EqR16z}, the formulas in this section are conjectural.
Again using \eqn{EqBy}, we find that
\beql{EqR16z}
\rho^{(m)} (n) =
\frac{m(n+1+2u) + \beta^{(m+1)}(n+1)}{m+1}
\eeq
for $0 \le n \le (m+2)^2-1$, where $u = \lfloor n/(m+2) \rfloor$.
Eliminating $\beta^{(m+1)}(n+1)$ from \eqn{REC1} and \eqn{EqR16z}
we obtain
\beql{EqR16zz}
\rho^{(m)} (n+1) = (m+2)\rho^{(m)} (n) + \sigma^{(m+1)}(n+1)
-m(n+1+2u) \, .
\eeq

In particular, for $0 \le n \le m+1$, $u=0$ and so \eqn{EqR16z} gives
\beql{EqR16}
\rho^{(m)} (n) =
\frac{(m+2)^{n+1} + (n+1)m(m+1) -1}{(m+1)^2} \, ,
\eeq
while for $m+2 \le n \le 2m+3$, $u=1$ and
\beql{EqR16a}
\rho^{(m)} (n) =
\frac{(m+2)^{n+1} + 2(m+1)(m+2)^{n-m-2} + (n+3)m(m+1) -1}{(m+1)^2} \, .
\eeq
Equation \eqn{EqR16} matches the smoothed values
on or below the diagonal $n=m+1$ in Table~\ref{T5},
\eqn{EqR16a} matches the values in the ``wedge'' $m+2 \le n \le 2m+3$,
and
\eqn{EqR16z} matches the values in the whole region bounded by the
``parabola'' $n \le (m+2)^2-1$.
Equation \eqn{EqR16a} is in fact a good estimate of $\rho^{(m)} (n)$
for all $m$ and $n$.

The greatest differences between the exact values $\pi^{(m)}(n)$
and the smoothed values $\rho^{(m)} (n)$ occur in the first row of 
Table~\ref{T5}.  
The ratio of terms $\pi^{(m)}(n+1)/\pi^{(m)}(n)$ in row $m$
of that table rapidly approaches $m+2$, and for fixed $m$
we find that
\beql{EqR13}
\pi^{(m)} (n) \approx \lambda_m (m+2)^n \,,
\eeq
where approximate values of $\lambda_m$ are
$$
\begin{array}{crcccccccc}
m &: & 1 & 2 & 3 & 4 & 5 & 6 & 7 & \cdots \\ 
\lambda_m & : & .778 & .447 & .312 & .240 & .194 & .163 & .140 & \cdots
\end{array}
$$
Curve-fitting suggests that
$$\lambda_m \approx
\frac{.956 m+2.11}{(m+1)^2}$$
which we approximate by
\beql{EqR14a}
\lambda_m \approx \frac{m+2}{(m+1)^2} \,,
\eeq
leading to
\beql{EqR14}
\pi^{(m)} (n) \approx \frac{(m+2)^{n+1}}{(m+1)^2}
\eeq
for $m$ fixed and $n$ large.

Since the leading terms in \eqn{EqR16z}, \eqn{EqR16} and \eqn{EqR16a}
agree with \eqn{EqR14},
we will take \eqn{EqR14} as our approximation 
to both $\pi^{(m)} (n)$ and $\rho^{(m)} (n)$
for all $m$ and $n$.

\subsection{ An estimate for $\tau^{(m)}(n)$ } \label{Sec42}
 From equations \eqn{defT}, \eqn{deftau}
we have
\beql{EqR18}
\tau^{(m)} (n+1) = \sum_{i=1}^n \sigma^{(m)} (i), \quad n \ge 1 \,.
\eeq

To simplify the analysis (we are only seeking a crude estimate in this section)
we suppose we have reached the end of block $B_n^{(m)}$ in $A^{(m)}$, where
$n= (m+1)^\mu + 1$ for some $\mu \ge 1$.
This block ends with the string $T_n^{(m)}$ of length $\tau^{(m)} (n)$.

Of the $(m+1)^\mu$ strings $S_i^{(m)}$, $1 \le i \le (m+1)^\mu$,
that appear in $T_n^{(m)}$, a fraction $\frac{m}{m+1}$ have $|i|_{m+1} =0$
and contribute
$\pi^{(m)} (0)$ to the sum;
a fraction $\frac{m}{(m+1)^2}$ have
$|i|_{m+1} =1$ and contribute $\pi^{(m)} (1)$;
and so on.
Therefore, from \eqn{EqR18},
$$\tau^{(m)} ((m+1)^\mu +1 ) =
(m+1)^\mu
\sum_{i=0}^{\mu -1}
\frac{m}{(m+1)^{i+1}} \pi^{(m)} (i) +\pi^{(m)} (\mu ) \,,$$
where the last term accounts for the final glue string
$S_{n-1}^{(m)}$.
Using \eqn{EqR14} this becomes
\beql{EqR20}
\tau^{(m)} ((m+1)^\mu +1) \approx 
\frac{m+2}{m+1} \,
\Big( (m+2)^\mu - m(m+1)^{\mu -1} \Big) \, .
\eeq

We summarize the discussion in the last two sections in the following:
\begin{obs}\label{Obs3}
After smoothing (cf. Observation \ref{Obs1}),
the sequence $\sm = \sm(1),
\sm(2), \ldots$ of the lengths of the glue strings in $A^{(m)}$ is a ruler
sequence given by \eqn{EqR12}, where the record values are given by
\eqn{EqR16z} (exact, for $n \le (m+2)^2-1$)
and \eqn{EqR14} (approximate, for all $m$ and $n$).
Equation \eqn{EqR20} gives an estimate for $\tau^{(m)}(n)$.
\end{obs}

\subsection{The first occurrence of $t$}\label{Sec43}
We use the above estimates to determine where a number $t \ge 1$
appears for the first time in $A$.
We already know from Tables~\ref{T1} and \ref{T2} that a 1 appears at
position 1, a 2 at position 3, a 3 at position 9, and a 4 at position 220, so we may
assume $t \ge 5$.

For fixed $t$, let $x(m)$ be the position where $t$ appears for 
the first time in the sequence $A^{(m)}$, for $1 \le m \le t$.
We will successively estimate $x(t)$, $x(t-1), \ldots, x(1)$, working backwards from
$$A^{(t)} = \underbrace{t,t, \ldots, t}_{t+1 ~\mbox{copies}}, t+1, \ldots, 
$$
where $t$ appears as the leading term,
and $t+1$ appears for the first time at position $t+2$.
A more detailed analysis of the beginning of $A^{(t)}$,
omitted here, shows that $t+2$ appears for the first time at position
\beql{EqS1}
\frac{(t+1)^{t+2} +2t-1}{t} \,.
\eeq
For example, 3 appears in $A$ at position 9, 4 in $A^{(2)}$ at position 42,
and 5 in $A^{(3)}$ at position 343.
Thus $x(t) =1$, $x(t-1) = t+1$, and
\beql{EqS2}
x(t-2) = \frac{(t-1)^t + 2t-5}{t-2} \,.
\eeq

We first consider the case $t=5$.
Since $x(3) = 343$, a 5 appears in $A^{(2)}$ for the first 
time at the end of block $B_i^{(2)}$, where $i$ is such that $\tau^{(2)} (i) =343$.
That is, $i$ is determined (see \eqn{EqR18} and Table~\ref{T4}) by the equation
$$
\tau^{(2)} (i) =
\underbrace{1+1+3+1+1+3+1+1+9+\cdots}_{i-1~\mbox{terms}} = 343 \,.
$$
By direct calculation, $i=80$, and again by direct calculation from \eqn{EqR17a},
a 5 appears in $A^{(2)}$ at position
$$
x(2) = \beta^{(2)} (80) =
77709404388415370160829246932345692180\,,
$$
which is $10^{37.9\ldots}$.
So 5 appears in $A = A^{(1)}$ at the end of block $B_i^{(1)}$, where $i$ is such that
$\tau^{(1)} (i) = x(2)$.
Setting $m=1$ in \eqn{EqR20} we get
$$
\tau^{(1)} (2^\mu +1 ) = \frac{3}{2} \cdot 3^\mu 
\Big( 1 - \frac{1}{2} \Big(\frac{2}{3}\Big)^\mu\Big) = x(2) \,,
$$
hence $\mu = 79.0\ldots$, $i=2^{79.0\ldots}$.
Setting $m=1$, $n = 2^{79.0\ldots}$ in 
\eqn{EqBc} we finally obtain
$$
x(1) = \beta^{(1)}(2^{79.0\ldots}) = \varepsilon _1 2^{2^{79.0\ldots}} = 10^{10^{23.3\ldots}} 
$$
for the position of the first $5$.

Consider now a general value of $t \ge 6$.
To find $x(t-3)$, we must solve (from \eqn{EqS2})
$$
\tau^{(t-3)} (i) = \frac{(t-1)^t + 2t-5}{t-2} \,.
$$
Setting $i= (t-2)^{\mu} +1$ and using \eqn{EqR20} we get
$$
\frac{t-1}{t-2} \,
\Big( (t-1)^\mu - (t-3)(t-2)^{\mu -1} \Big) =
\frac{(t-1)^t + 2t-5}{t-2} \,,
$$
hence $\mu \approx t-1$, and so
$$
x(t-3) = \beta^{(3)}
((t-2)^\mu +1) =
\varepsilon _3(t-2)^{(t-2)^{t-1}} \,.
$$
The $\varepsilon _3$ may be ignored, since it can be absorbed into the tower of
exponentials.
The next iteration gives
$$
x(t-4) = (t-3)^{(t-3)^{(t-2)^{t-1}}}$$
and eventually we obtain
\beql{EqAns}
x(1) = 
2^{2^{\Snj 3^{\Snj 4^{\cdot^{\cdot^{\cdot^{\Snj {t-1}}}}}}}} \,,
\eeq
a tower of height $t-1$.
We formalize this as our final
\begin{obs}\label{Obs4}
The number $t \ge 5$ appears for the first time
in the sequence $A$ at about position \eqn{EqAns}.
\end{obs}

\noindent{\bf Remark.}
The position of the first $5$ can be estimated more accurately
using the tabular construction of Section \ref{Sec4BB},
thus avoiding the simplifying assumptions made in Section \ref{Sec42}.
This analysis predicts that the first $5$
will appear at about position
$$
\varepsilon _1 \times 2^{418090195952691922788353} = 10^{10^{23.09987\ldots}} \,. 
$$
where $\varepsilon _1$ is given in \S\ref{Sec4B}.
We omit the details.

\section{Comments and generalizations}\label{S5}
\subsection{The Finiteness Conjecture}\label{Sec51}
The proof of Theorem \ref{Th1} would have been simpler if we 
had known in advance that the glue strings $S_n^{(m)}$ were finite.
This would follow from the following:

\paragraph{Finiteness Conjecture.}
{\em For integers $m \ge 2$ and $r \ge 1$, let $x(1), x(2),\ldots, x(r)$ be a
string from $\PP_m^+$.
Let $x(n+1) = \sC (x(1), x(2), \ldots, x(n))$ for $n \ge r$.
Then for some $n \ge r+1$, $x(n) < m$.} \medskip

In other words, there is no finite starting string from $\PP_m^r$ which 
extends under repeated application of the map $\sC$ to an infinite sequence from
$\PP_m^\ast$.
Sooner or later a term less than $m$ must appear.

Although this conjecture seems very plausible, we have not been able to find a proof.
If one tries to construct a starting string which extends for a 
long time without dropping below $m$ one quickly runs into difficulties.
Let $m_1$ and $m_2$ be respectively the smallest and largest values in 
the starting string.
Then no number outside the range $[m_1, m_2]$ --- or in fact 
any number not in the starting string --- may appear in the 
resulting sequence,
for such a number is immediately followed by a 1, terminating the sequence.
So if the sequence is infinite it must be bounded.

As an experiment we considered all $2^n$ starting strings of length $n \le 30$
consisting just of 2's and 3's, and calculated the lengths of the 
resulting strings until just before the appearance of the first 1.
The maximum and average lengths are shown in Table~\ref{T8}.
The average length seems to approach $n \, +$ constant,
but the maximum length is harder to understand, and it would be nice to have more data.
Does the maximum length continue to grow linearly,
or are there further jumps of ever-increasing size?
We do not know.
\begin{table}[htb]
$$
\begin{array}{|r|r|r@{}l|} \hline
n & \mbox{Max.} && \mbox{Ave.} \\ \hline
1 & 1 & 1& \\
2 &  4 & 2&.75 \\
3 &  5 & 3&.75 \\
4 &  8 & 5&.125 \\
5 &  9 & 6&.2187 \\
6 & 14 & 7&.5 \\
7 & 15 & 8&.5703 \\
8 & 66 & 10&.2734 \\
9 & 68 & 11&.3828 \\
10 & 70 & 12&.5293 \\
11 & 123 & 13&.6099 \\
12 & 124 & 14&.6658 \\
13 & 125 & 15&.6683 \\
14 & 132 & 16&.6957 \\
15 & 133 & 17&.7047 \\
16 & 134 & 18&.7168 \\
17 & 135 & 19&.7206 \\
18 & 136 & 20&.7278 \\
19 & 138 & 21&.7304 \\
20 & 139 & 22&.7341 \\
21 & 140 & 23&.7353 \\
22 & 142 & 24&.7372 \\
23 & 143 & 25&.7379 \\
24 & 144 & 26&.7388 \\
25 & 145 & 27&.7391 \\
26 & 146 & 28&.7396 \\
27 & 147 & 29.&7398 \\
28 & 148 & 30.&74 \\
29 & 149 & 31.&7401 \\
30 & 150 & 32.&7402 \\
\hline
\end{array}
$$
\caption{Maximum and average length of string produced by any starting 
sequence of $n$ 2's and 3's, stopping when first 1 is reached.}
\label{T8}
\end{table}

Table~\ref{T9} shows the starting strings of lengths 2, 4, 6, 8 and 11
(when there are jumps in the maximum length)
and the strings of record lengths 4, 8, 14, 66 and 123 that they produce.
These five starting strings are unique.
\begin{table}[htb]
$$
\begin{array}{ll}
n = 2 & \mbox{Starting string 2 2 produces length 4:} \\
      & \mbox{2 2 2 3} \\ [+.2in]
n = 4 & \mbox{Starting string 2 3 2 3 produces length 8:} \\
      & \mbox{2 3 2 3 2 2 2 3} \\ [+.2in]
n = 6 & \mbox{Starting string 2 2 2 3 2 2 produces length 14:} \\
      & \mbox{2 2 2 3 2 2 2 3 2 2 2 3 3 2} \\[+.2in]
n = 8 & \mbox{Starting string 2 3 2 2 2 3 2 3 produces length 66:} \\
      & \mbox{2 3 2 2 2 3 2 3 2 2 2 3 2 2 2 3 2 2 3 2 2 2 3 2 2 2 3 2 3 2} \\
      & \mbox{2 2 3 2 2 2 3 2 2 3 2 2 2 3 2 2 2 3 2 3 2 2 2 3 2 2 2 3 2 2} \\
      & \mbox{3 2 2 3 3 2} \\ [+.2in]
n = 11 & \mbox{Starting string 2 2 3 2 3 2 2 2 3 2 2 produces length 123:} \\
       & \mbox{2 2 3 2 3 2 2 2 3 2 2 2 3 2 2 3 2 2 2 3 2 2 2 3 2 3 2 2 2 3} \\
       & \mbox{2 2 2 3 2 2 3 2 2 2 3 2 2 2 3 2 3 2 2 2 3 2 2 2 3 2 2 3 2 2} \\
       & \mbox{2 3 2 3 2 2 2 3 2 2 2 3 2 2 3 2 2 2 3 2 2 2 3 2 3 2 2 2 3 2} \\
       & \mbox{2 2 3 2 2 3 2 2 2 3 2 2 2 3 2 3 2 2 2 3 2 2 2 3 2 2 3 2 3 2 2 2 3}
\end{array}
$$
\caption{Starting strings of $n$ 2's and 3's which extend under 
the map $\sC$ for a record number of steps before reaching a 1.
Each entry also shows the final string (until just before the first $1$
is reached).}
\label{T9}
\end{table}

\subsection{Curling number transforms of other sequences}\label{Sec52}

It is interesting to apply the curling number transform to other sequences,
particularly those for which the definition involves properties of substrings.
For example, the binary Thue-Morse sequence
$$
\begin{array}{ccccccccccccccccc}
0 & 1 & 1 & 0 & 1 & 0 & 0 & 1 & 1 & 0 & 0 & 1 & 0 & 1 & 1 & 0 \\ 
1 & 0 & 0 & 1 & 0 & 1 & 1 & 0 & 0 & 1 & 1 & 0 & 1 & 0 & 0 & 1 & \ldots
\end{array}
$$
(A10060 in \cite{OEIS}) has the property that it contains no cubes $UUU$ as substrings
(see \cite{AlSh03}, \cite{Lot83}, \cite{Sal81}, \cite{OEIS} for further
information).
Its curling number transform, which naturally contains only 1's
and 2's, is
$$
\begin{array}{ccccccccccccccccc}
1 & 1 & 1 & 2 & 1 & 1 & 2 & 2 & 1 & 2 & 1 & 2 & 2 & 1 & 2 & 2 \\ 
1 & 2 & 2 & 2 & 1 & 2 & 2 & 2 & 2 & 2 & 1 & 2 & 2 & 1 & 2 & 2 & \ldots 
\end{array}
$$
(A93914).
We leave it to the interested reader to investigate the properties of
this sequences and of the other new sequences mentioned below.

There are many examples of ternary sequences which contain no squares,
and of course their curling number transforms are simply the all-ones sequence {\bf 1}.
However, the lexicographically earliest sequence from $\PP^\ast$ whose 
transform is {\bf 1} is the ruler sequence $r$ (A1511)
mentioned in Section \ref{Sec41}.

We give one further example.
The Kolakoski sequence is a sequences of $1$'s and $2$'s
defined by $K(1) =1$, $K(n) =$ length of $n$-th run:
$$
\begin{array}{ccccccccccccccccc}
1 & 2 & 2 & 1 & 1 & 2 & 1 & 2 & 2 & 1 & 2 & 2 & 1 & 1 & 2 & 1 \\ 
1 & 2 & 2 & 1 & 2 & 1 & 1 & 2 & 1 & 2 & 2 & 1 & 1 & 2 & 1 & 1 &\ldots
\end{array}
$$
(A2 in \cite{OEIS}).
This also contains no cubes.
The transformed sequence is
$$
\begin{array}{ccccccccccccccccc}
1 & 1 & 1 & 2 & 1 & 2 & 1 & 1 & 2 & 2 & 1 & 2 & 2 & 2 & 2 & 1  \\ 
1 & 2 & 2 & 2 & 1 & 1 & 2 & 2 & 1 & 2 & 2 & 2 & 1 & 2 & 1 & 1 & \ldots
\end{array}
$$
(A93921).

\subsection{Generalizations}\label{Sec53}
In this final section we briefly mention a few of the possible 
generalizations of the sequence $A$. 
\medskip

\noindent
(i)~The recurrence \eqn{EqM2} may be replaced by $a(1) =1$, $a(n+1) =
f(\sC (a(1), \ldots, a(n)))$ for $n \ge 1$, for any suitable function $f$.
For example, $f(x) = \mbox{floor} (x/2 )$ produces
$$
\begin{array}{ccccccccccccccccc}
0 & 0 & 1 & 0 & 0 & 1 & 1 & 1 & 1 & 2 & 0 & 0 & 1 & 0 & 0 & 1 \\ 
1 & 1 & 1 & 2 & 1 & 0 & 0 & 1 & 0 & 0 & 1 & 1 & 1 & 1 & 2 & 0 &\ldots
\end{array}
$$
(A91970), which presumably has an even slower rate of growth than $A$.
\medskip

\noindent
(ii)~$a(1) = a(2) =1$,
$a(n+2) = \sC (a(1),\ldots, a(n))$ for $n \ge 1$ produces
$$
\begin{array}{ccccccccccccccccc}
1 & 1 & 1 & 2 & 3 & 1 & 1 & 1 & 2 & 3 & 1 & 2 & 2 & 1 & 2 & 1 \\ 
1 & 2 & 2 & 1 & 2 & 1 & 1 & 2 & 2 & 2 & 2 & 3 & 4 & 1 & 1 & 1 &\ldots
\end{array}
$$
(A94006).
This has the property that its curling number transform is the same 
sequence but shifted one place to the left.
\medskip

\noindent
(iii)~A greedy version of $A^{(2)}$. Let $g(1) = 2$. For $n \ge 1$, let
$k = \sC (g(1), \ldots , g(n))$. If $k > 1$, $g(n+1) = k$ (as in $A^{(2)}$),
but if $k=1$, choose $g(n+1)$ so that $\sC (g(1), \ldots , g(n+1))$ is
maximized. If there is more than one choice for $g(n+1)$, pick
the smallest. The resulting sequence (A94321) is:
$$
\begin{array}{ccccccccccccccccc}
2 & 2 & 2 & 3 & 3 & 2 & 2 & 2 & 3 & 3 & 2 & 2 & 2 & 3 & 2 & 2 \\
2 & 3 & 2 & 2 & 2 & 3 & 3 & 2 & 2 & 2 & 3 & 2 & 2 & 2 & 3 & 2 & \ldots
\end{array}
$$
\medskip

\noindent
(iv)~A two-dimensional version of $A$. 
Define $t(i,j)$, $i \ge 1$, $j \ge 1$, as follows:
$t(i,1) = t(1,i) = a(i)$. For $i, j \ge 1$,
$t(i+1,j+1) = \max \{ k_1, k_2 \}$,
where $k_1 = \sC (t(i+1,1), t(i+1,2), \ldots, t(i+1,j))$,
$k_2 = \sC (t(1,j+1), t(2,j+1), \ldots, t(i,j+1))$
--- see Table~\ref{TabTT} (A94781). 
The first two rows (or columns) give $A$ and the third row (or column) is $A^{(2)}$.

\begin{table}[htb]
$$
\begin{array}{cccccccccccccccc}
1 & 1 & 2 & 1 & 1 & 2 & 2 & 2 & 3 & 1 & 1 & 2 & 1 & 1 & \ldots \\
1 & 1 & 2 & 1 & 1 & 2 & 2 & 2 & 3 & 1 & 1 & 2 & 1 & 1 & \ldots \\
2 & 2 & 2 & 3 & 2 & 2 & 2 & 3 & 2 & 2 & 2 & 3 & 3 & 2 & \ldots \\
1 & 1 & 3 & 1 & 1 & 3 & 3 & 2 & 1 & 1 & 2 & 1 & 1 & 2 & \ldots \\
1 & 1 & 2 & 1 & 1 & 2 & 2 & 2 & 3 & 1 & 2 & 1 & 1 & 2 & \ldots \\
2 & 2 & 2 & 3 & 2 & 1 & 1 & 2 & 1 & 2 & 3 & 2 & 2 & 3 & \ldots \\
. & . & . & . & . & . & . & . & . & . & . & . & . & . & \ldots \\
\end{array}
$$

\caption{A two-dimensional version of the sequence.}
\label{TabTT}
\end{table}
\medskip

\noindent
(v)~J. Taylor \cite{Tay04} has suggested two broad generalizations 
of the original recurrence.
Let $\sim$ be an equivalence relation on strings of integers of each fixed length.
Write
\beql{EqT1}
a(1) a(2) \ldots a(n) = XY_1 Y_2 \ldots Y_k \,,
\eeq
where the $Y_i$ are nonempty strings with $Y_1 \sim Y_2 \sim \cdots \sim Y_k$
and $k$ is maximal; then $a(n+1) =k$.
Choosing $\sim$ to be the identity relation gives $A$.
Taylor has contributed several interesting generalizations of $A$ to \cite{OEIS} 
obtained from other equivalence relations.
For example, if two strings are equivalent if one is a permutation
of the other, the resulting sequence is
$$
\begin{array}{ccccccccccccccccc}
1 & 1 & 2 & 1 & 1 & 2 & 2 & 2 & 3 & 1 & 1 & 2 & 1 & 1 & 2 & 2 \\ 
2 & 3 & 2 & 2 & 2 & 3 & 2 & 2 & 2 & 3 & 3 & 2 & 2 & 4 & 1 & 1 & \ldots
\end{array}
$$
(A91976), which agrees with $A$ for the first $19$ terms.
But after
$$
\begin{array}{ccccccccccccccccccc}
1 & 1 & 2 & 1 & 1 & 2 & 2 & 2 & 3 & 1 & 1 & 2 & 1 & 1 & 2 & 2 & 2 & 3 & 2
\end{array}
$$
the next term is now 2, not 1, since we can take $U=1~1~2~1~1$,
$Y_1 = 2~2~2~3~1~1~2$,
$Y_2= 1~1~2~2~2~3~2$,
where $Y_2$ is a permutation of $Y_1$.
 
\noindent
(vii)~More generally, Taylor suggests using a partial order $\prec$ on 
integer strings of all lengths (not just strings of the same length),
and requiring the $Y_i$
in \eqn{EqT1} to be nonempty and satisfy
$Y_1 \prec Y_2 \prec \cdots \prec Y_k$ where $k$ is maximal.
For further examples of Taylor's sequences the reader is referred to the 
entries A91975 and A92331--A92335 in \cite{OEIS}.

\subsection*{Acknowledgements}
We thank J. Taylor for telling us about his generalizations of the sequence.

\end{document}